\documentclass[11pt]{article}

\usepackage{epsfig}
\usepackage{graphicx}
\usepackage[latin1]{inputenc}

\usepackage{amsfonts,amsthm,bbm,amssymb,epsf,amsmath,%
latexsym,amscd,amsfonts,enumerate,amsthm,supertabular,color,url}

\usepackage{subfigure,psfrag}

\newcommand{\p}{\partial}

\usepackage{algorithm}
\usepackage{algorithmic}

\newcommand{\mpr}[1]{\relax}
\newcommand{\mprup}[1]{\relax}
\newcommand{\mprdown}[1]{\relax}

\newtheorem{theorem}{Theorem}[section]
\newtheorem{lemma}[theorem]{Lemma}
\newtheorem{corollary}[theorem]{Corollary}
\newtheorem{proposition}[theorem]{Proposition}
\newtheorem{definition}[theorem]{Definition}

\newcommand{\pf}{{\noindent \bf Proof.} }
\newcommand{\eop}{{ \vrule height7pt width7pt depth0pt}\par\bigskip}

\def\sp{\,:\;}
\def\ie{{\sl \thinspace i.e.},\ }

\def\tri{\triangle}

\newcommand{\RR}{\mathbb{R}}

\def\uw{\tilde u}

\DeclareMathOperator{\supp}{supp}

\DeclareMathOperator{\st}{star}
\DeclareMathOperator{\diag}{diag}

\begin{document}

\title{Polynomial Finite Element Method for Domains Enclosed by Piecewise
Conics}%

\author{Oleg Davydov\thanks{Department of
Mathematics, University of Giessen,
    Department of Mathematics,
    Arndtstrasse 2, 
    35392 Giessen,
    Germany, {\tt oleg.davydov@math.uni-giessen.de}}
\and Georgii Kostin\thanks{Institute for Problems in Mechanics,
Russian Academy of Sciences, 101-1, Vernadskogo pr., Moscow, 117526, Russia,
{\tt kostin@ipmnet.ru}}
\and Abid Saeed\thanks{Department of Mathematics,
Kohat University of Science and Technology,
Kohat, Pakistan, {\tt abidsaeed@kust.edu.pk}}}

\date{August 21, 2015}

\maketitle

\begin{abstract}
We consider bivariate piecewise polynomial finite element spaces  for 
curved domains bounded by piecewise conics satisfying homogeneous boundary conditions,
construct stable local bases for them using Bernstein-B\'ezier techniques, prove error bounds and
develop optimal assembly algorithms for the finite element system matrices.
Numerical experiments confirm the effectiveness of the method.
\end{abstract}

\section{Introduction\label{intro}}\mprup{intro}

Spaces of multivariate piecewise polynomial splines are usually 
defined on triangulated polyhedral domains without imposing any boundary
conditions. However, applications such as the finite element method require at
least the ability to prescribe zero values on parts of the boundary. Fitting data
with curved discontinuities of the derivatives is another
situation where the interpolation of prescribed values along a lower dimensional manifold
is highly desirable. It turns out that such conditions make the
otherwise well understood spaces of e.g.\ bivariate $C^1$ macro-elements on triangulations 
significantly more complex. Even in the simplest case of a polygonal domain, the
dimension of the space of splines vanishing on the boundary is dependent on its geometry, with 
consequences for  the construction of stable bases (or stable minimal determining sets) 
\cite{DSa12,DSa13}.

Since splines are piecewise polynomials, it is convenient to model curved features by  
piecewise algebraic surfaces so that the spline space naturally splits out the subspace of functions
vanishing on such a surface. Indeed, implicit algebraic surfaces are a well-established modeling tool
in CAGD \cite{ImplSurf}, and the ability to exactly reproduce some of them (e.g.\ circles 
or cylinders) is a highly desirable feature for any modeling method \cite{FarinBook}. 

On the other hand, the finite element analysis benefits a lot from the isogeometric approach
\cite{HCB05}, where the geometric models of the boundary are used exactly in the
form they exist in a CAD system rather than undergoing a remeshing to fit into the traditional
isoparametric finite element scheme. While the isogeometric analysis introduced in \cite{HCB05} is
based on the most widespread modeling tool of NURBS and benefits from the many attractive features of
tensor-product B-splines, it also inherits some of their drawbacks, such as complicated local 
refinement (see for example \cite{BCS10}). 

In this paper we explore an isogeometric method which combines modeling with algebraic curves with the
standard triangular piecewise polynomial finite elements in the simplest case of planar domains defined by
piecewise quadratic algebraic curves (conic sections). 
Remarkably, the standard Bernstein-B\'ezier  
techniques for dealing with piecewise polynomials on triangulations \cite{LSbook,Sch15} as well as recent
optimal assembly algorithms \cite{AAD11,AAD15,ADS} for high order elements can be carried over to this case without significant loss of
efficiency. Some of the material, especially in Sections 4 and 6 is based on the thesis \cite{Abid_thesis} 
of one of the authors. Note that we only consider $C^0$ elements for elliptic problems with homogenous Dirichlet boundary 
conditions, although preliminary results on a direct implementation of non-homogenous 
Dirichlet boundary conditions can be found in \cite{Abid_thesis}. %

In contrast to both the isoparametric curved finite elements and the isogeometric analysis, our approach 
does not require parametric patching on curved subtriangles, and therefore does not depend on the invertibility
of the Jacobian matrices of the nonlinear geometry mappings. Therefore our finite elements remain piecewise polynomial
everywhere in the physical domain. This in particular facilitates a relatively straightforward
extention  to $C^1$ elements on piecewise conic domains, which have also been considered in \cite{Abid_thesis} and tested numerically on the approximate solution  
of fully nonlinear elliptic equations by Böhmer's method \cite{Boehmer08}. Full details of the theory of 
these elements are postponed to our forthcoming paper \cite{DSaC1}. 

There are some connections to the  weighted extended B-spline (web-spline) method  \cite{HRW01}. In particular, in our error
analysis we use a technical lemma (Lemma~\ref{Hardyl}) proved in  \cite{HRW01}. Indeed, the quadratic
polynomials that define the curved edges of the pie-shaped triangles at the domain boundary
are factored out  from the local polynomial spaces and hence act as weight functions on certain subdomains.
They remain however integral parts of the spline spaces in our case and are generated naturally from the conic
sections defining the domain, thus bypassing the problem of the computation of a smooth global  weight function
needed in the web-spline method.

The paper is organized as follows. We introduce in Section 2 the spaces $S_{d,0}(\tri)$ of $C^0$ piecewise
polynomials of degree $d$ on
domains bounded by a number of conic sections, with homogeneous boundary conditions and investigate in Section 3 their
approximation power for functions in Sobolev spaces $H^m(\Omega)$ vanishing on the boundary, which
leads in particular to the error bounds in the form $\mathcal{O}(h^{m})$ in the $L_2$-norm and 
$\mathcal{O}(h^{m-1})$ in the $H^1$-norm for the solutions of elliptic problems by the Ritz-Galerkin finite element
method. Section 4 is devoted to the development of a basis for $S_{d,0}(\tri)$ of Bernstein-B\'ezier type
important for a numerically stable and  efficient implementation of the method. Some implementation issues
specific for the curved elements are treated in Section 5, including the fast assembly of the system matrices.
Finally, Section 6 presents several numerical experiments involving the Poisson problem 
on two different curved domains, as well as the circular membrane eigenvalue problem. The results confirm 
the effectiveness of our method both in $h$- and $p$-refinement settings.

\section{Piecewise polynomials on  piecewise conic domains\label{spacesC0}}\mprup{spaces}

Let $\Omega\subset\RR^2$ be a bounded curvilinear polygonal domain with 
$\Gamma=\partial \Omega=\bigcup_{j=1}^n\overline{\Gamma}_j$, 
where each $\Gamma_j$ is an open arc of an algebraic curve of
at most second order
(\ie either a straight line or a conic). For simplicity we assume that 
$\Omega$ is simply connected.
Let $Z=\{z_1,\ldots,z_n\}$ be the set of the endpoints of all
arcs numbered counter-clockwise such that $z_j,z_{j+1}$ are the
endpoints of $\Gamma_j$, $j=1,\ldots,n$. (We set $z_{j+n}=z_j$.) 
Furthermore, for each $j$ we denote by $\omega_j$ the internal 
angle between the tangents $\tau^+_j$ and $\tau^-_j$ to $\Gamma_j$
and $\Gamma_{j-1}$, respectively, at $z_j$. We assume that
$0<\omega_j\le 2\pi$,  %
and set $\omega:=\min\{\omega_j\sp 1\le j\le n\}$.

Our goal 
is to develop an $H^1$-conforming finite element method 
with polynomial shape functions 
suitable for solving second order elliptic problems on curvilinear polygons
of the above type. 

\begin{figure}
\centering
\begin{tabular}{c}
\includegraphics[height=0.3\textwidth]{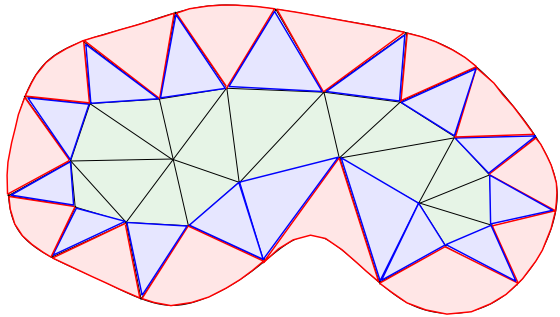} 
\end{tabular}
\caption{
A triangulation of a curved domain with ordinary triangles (green), 
pie-shaped triangles (pink) and buffer triangles (blue).
}\label{curved_Alltriangles}
\end{figure}

Let $\tri$ be a \emph{triangulation} of $\Omega$, \ie a 
subdivision of $\Omega$ into triangles, where 
each triangle $T\in\tri$ has at most one edge replaced with a curved segment of 
the boundary $\partial\Omega$, and the intersection of any pair of the triangles
is either a common vertex or a common (straight) edge if it is non-empty.
The triangles with a curved edge are said to be \emph{pie-shaped}.
Any triangle $T\in\tri$ that shares at least one edge with a pie-shaped triangle
is called a \emph{buffer} triangle, and the remaining triangles are
\emph{ordinary}. We denote by $\tri_0$, $\tri_B$ and $\tri_P$ the sets of all
ordinary, buffer and pie-shaped triangles of $\tri$, respectively. Thus,
$$
\tri=\tri_0\cup\tri_{B} \cup \tri_P$$
is a disjoint union, see Figure~\ref{curved_Alltriangles}.
We emphasize that a triangle with only straight edges on the boundary of
$\Omega$ does not belong to $\tri_P$.

We denote by $\mathbb{P}_d$  the space of all bivariate polynomials of total degree
at most $d$. For each  $j=1,\ldots,n$, let $q_j\in\mathbb{P}_2$ be a polynomial such that 
$\Gamma_j\subset\{x\in\RR^2\sp q_j(x)=0\}$. By multiplying $q_j$ by $-1$ if needed, we ensure that
 $\p_{\nu_x} q_j(x)<0$ for all $x$ in the interior of
$\Gamma_j$, where $\nu_x$ denotes the unit outer normal to the boundary at $x$,
and $\p_{a}:=a\cdot\nabla$ is the directional derivative with respect to a vector $a$. Hence,
$q_j(x)$ is positive for points in $\Omega$ near the boundary segment $\Gamma_j$.
We assume that $q_j\in \mathbb{P}_1$ or $q_j\in \mathbb{P}_2\setminus\mathbb{P}_1 $ 
depending on whether $\Gamma_j$ is a straight interval or a genuine conic arc.

Furthermore, let $V,E,V_I,E_I,V_B$ and $E_B$ denote the set of all vertices, 
all edges, interior vertices, interior edges, boundary vertices and boundary
edges of $\tri$, respectively. 
For each $v\in V$, $\st(v)$ stands for the union of all triangles in $\tri$
attached to $v$. We also denote by $\theta$ the smallest angle
of the triangles $T\in\tri$, where the angle between an interior
edge and a boundary segment is understood in the tangential sense.

We assume that $\tri$ satisfies the following \emph{Conditions}:
\begin{itemize}
\item[(a)] $Z=\{z_1,\ldots,z_n\}\subset V_B$.
\item[(b)] No interior edge has both endpoints on the boundary.
\item[(c)] %
If $q_{j}/q_{j-1}\ne \hbox{const}$ and at least one of 
$q_{j},q_{j-1}$ belongs to $\mathbb{P}_2 \backslash \mathbb{P}_1 $,
then there is at least one triangle $T\in\tri_{B}$ attached to $z_j$.
\item[(d)] Every $T\in\tri_P$ is star-shaped with respect to its interior vertex $v$.
\item[(e)]  For any  $T\in\tri_P$ with its curved side on $\Gamma_j$,
$q_j(z)>0$ for all $z\in T\setminus \Gamma_j$.
\end{itemize}
Note that (b) and (c) can always be achieved by a slight modification of
a given triangulation, while (d) and (e) hold for sufficiently fine triangulations. 

\begin{figure}
\centering
\begin{tabular}{c}
\psfrag{R}{}
\psfrag{A}{}
\includegraphics[height=0.35\textwidth]{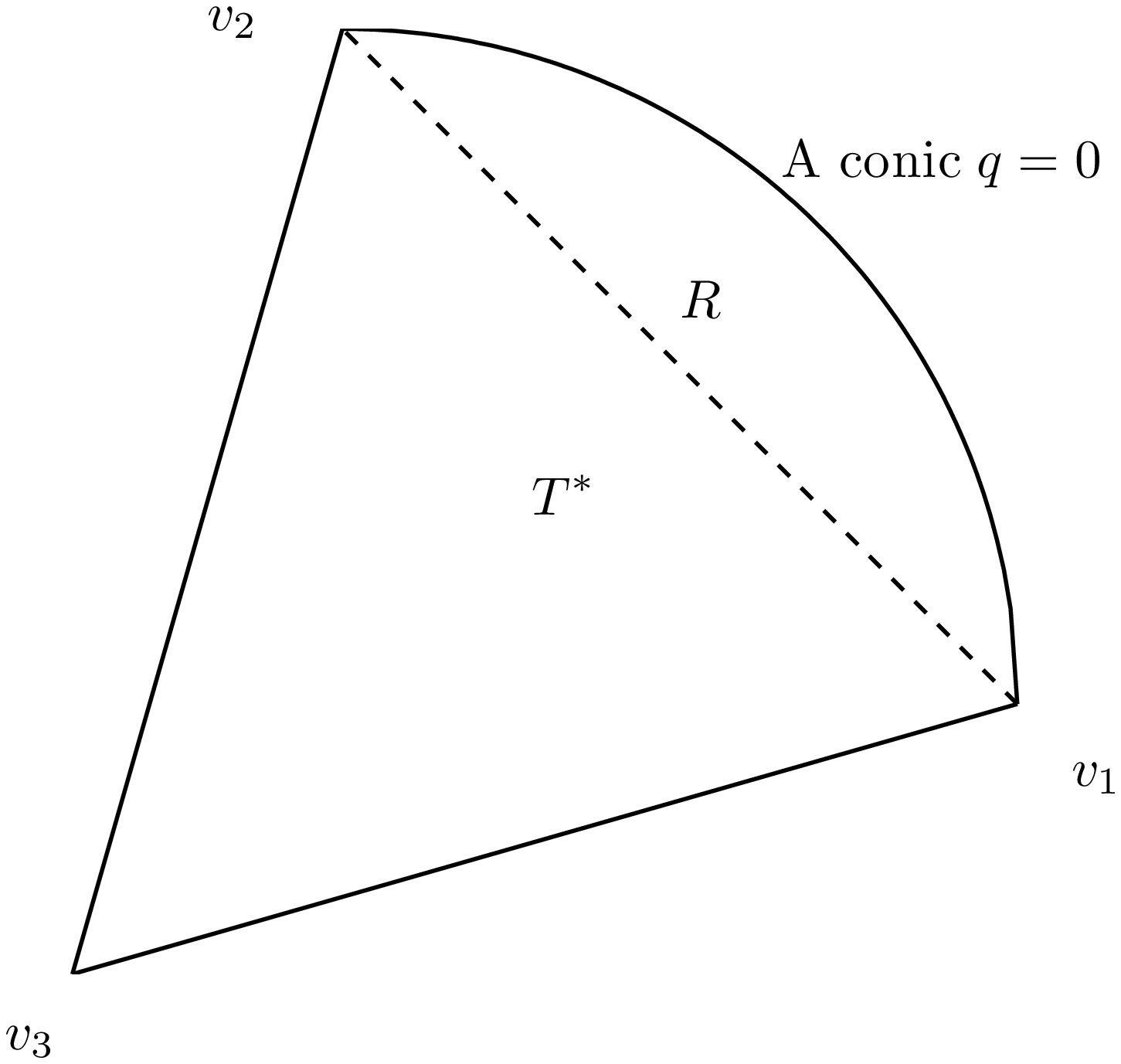} 
\psfrag{R}{}
\psfrag{A}{}
\includegraphics[height=0.35\textwidth]{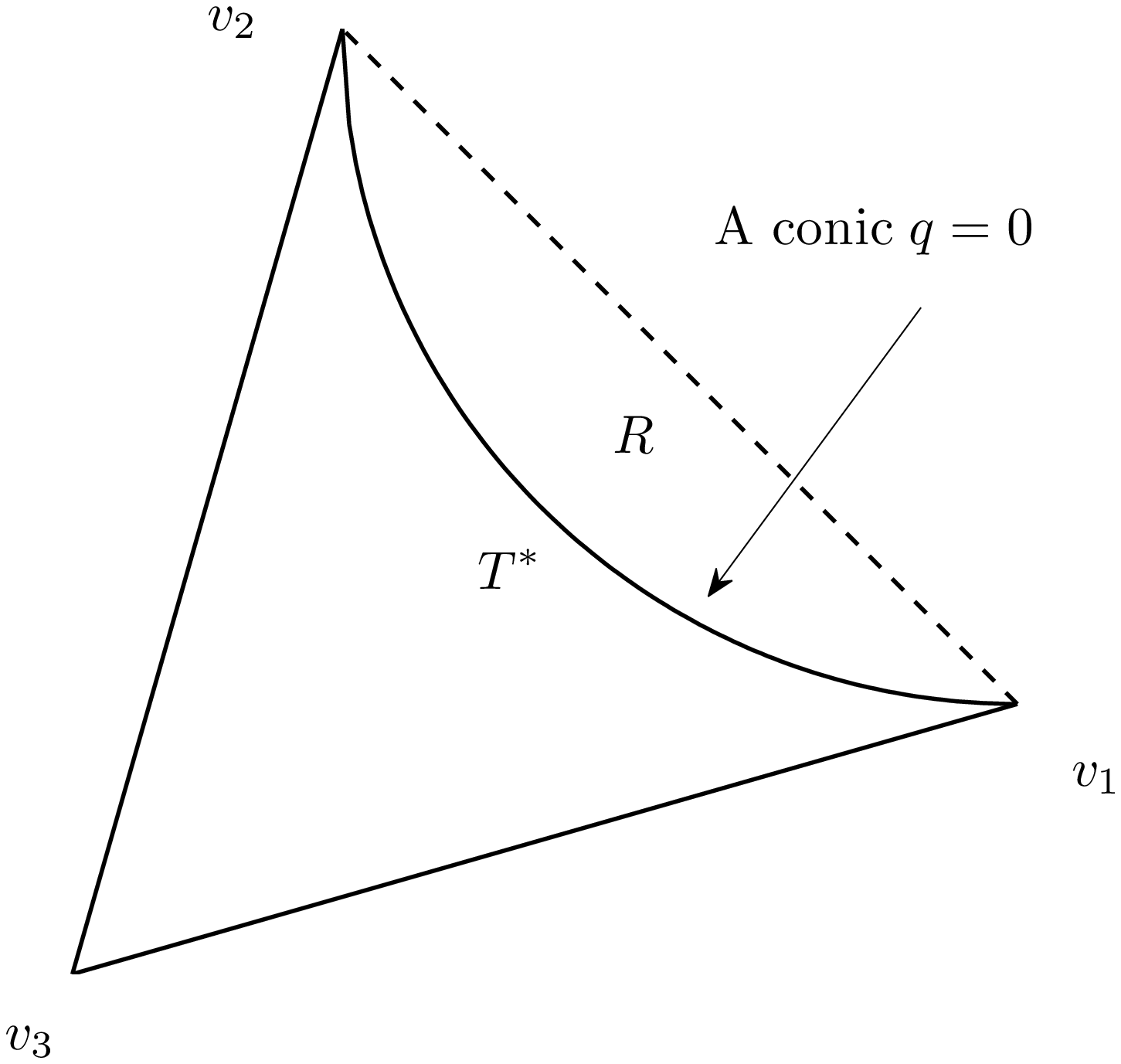} 
\end{tabular}
\caption{
A pie-shaped triangle with a curved edge and the
associated triangle $T^\ast$ with straight sides and vertices $v_1,v_2,v_3$.
The curved edge can be either outside (left) or inside $T^\ast$ (right).
}\label{pieshapedT1}
\end{figure}

For any $T\in\tri$, let $h_T$ denote the diameter of $T$, and let
$\rho_T$ be the radius of the disk $B_T$ inscribed in $T$ if $T\in\tri_0\cup\tri_B$ or in $T\cap T^*$ if 
$T\in\tri_P$, where $T^\ast$ denotes the triangle obtained by joining the boundary vertices of $T$
by a straight line, see Figure~\ref{pieshapedT1}. Note that every triangle $T\in\tri$ is star-shaped
with respect to $B_T$ in the sense of \cite[Definition 4.2.2]{BrennerScott}. In particular, for
$T\in\tri_P$ this follows from Condition (d) and the fact that the  conics do not possess inflection points.

We assume that $R,A,B$ are positive constants such that 
\begin{align}
&h_T\le R\rho_T,\quad \forall T\in \tri,\label{shape}
\end{align}
and, for any $T\in\tri_P$,
\begin{align}
& q_j(z)\leq Aq_j(v) ,\quad \forall z\in T,\label{q3}\\
&\p_{v-z}q_j(z)\geq Bq_j(v) ,\quad \forall z\in  T\cap\Gamma_j,\label{q4}
\end{align}
where $\Gamma_j$ is the conic arc containing the curved edge of $T$, and $v$ is the interior vertex of  $T$.
The constants $R,A,B$ exist for any triangulation $\tri$ of the above type and 
are responsible for the shape regularity of its triangles.

Let $d\ge1$. We set
\begin{align*}
S_d(\tri)&:=\{s\in C^0(\Omega)\sp s|_T\in
\mathbb{P}_{d+i},\;T\in\tri_i,\;i=0,1\},
\quad \tri_1:=\tri_P\cup\tri_B,
				 ,\\
S_{d,0}(\tri)&:=\{s\in S_d(\tri)\sp s|_\Gamma = 0\}. %
\end{align*}

\section{Error bounds \label{bounds}}\mprup{bounds}

In this section we provide some typical  error bounds for the spaces $S_{d,0}(\tri)$ in the context of 
the approximation theory and the finite element method. 

We denote by $\p^\alpha f$, $\alpha\in\mathbb{Z}^2_+$, the partial derivatives of $f$ and consider the
usual Sobolev spaces $H^m(\Omega)$ with the seminorm and norm defined by
$$
|f|_{H^m(\Omega)}^2=\sum_{|\alpha|=m}\|\p^\alpha f\|_{L^2(\Omega)}^2,\quad
\|f\|_{H^m(\Omega)}^2=\sum_{k=0}^m|f|_{H^k(\Omega)}^2\quad (H^0(\Omega)=L^2(\Omega)),$$
where $|\alpha|:=\alpha_1+\alpha_2$. We set
$H^1_0(\Omega)=\{f\in H^1(\Omega):f|_{\p\Omega}=0 \}$.

\subsection{Approximation error}
Let $J_1$ be the set of all  $j=1,\ldots,n$ such that $q_j\in\mathbb{P}_1$, and 
$J_2:=\{1,\ldots,n\}\setminus J_1$. 
For each $j\in J_2$ we choose a domain $\Omega_j\subset\Omega$ with Lipschitz boundary such that 
\begin{itemize}
\item[(a)] $\p\Omega_j\cap\p\Omega = \Gamma_j$,
\item[(b)] $\p\Omega_j\setminus\p\Omega$ is composed of a finite number of straight line segments,
\item[(c)]  $q_j(x)>0$ for all $x\in \overline{\Omega}_j\setminus \Gamma_j$, and
\item[(d)]  $\Omega_j\cap\Omega_k=\emptyset$ for all $j,k\in J_2$. 
\end{itemize}
We assume that the triangulation $\tri$ is such that for each $j\in J_2$,
\begin{equation} \label{q2}
\text{$\overline{\Omega}_j$ contains every triangle $T\in\tri_P$ whose curved edge is part of
$\Gamma_j$.}
\end{equation}
Note that \eqref{q2} will hold with the same set $\{\Omega_j: j\in J_2\}$
for any triangulations obtained by subdividing the triangles 
of $\tri$. In addition we assume in this section for the sake of simplicity of the analysis that
\begin{equation} \label{q1}
\text{no pair of pie-shaped triangles shares an edge.}
\end{equation}
This can be achieved e.g.\ by applying Clough-Tocher splits
\cite[Section 6.2]{LSbook} to certain pie-shaped triangles. Finally, without loss of generality we assume
that 
\begin{equation} \label{q0}
\max_{x\in\overline{\Omega}_j}\|\nabla q_j(x)\|_2\le 1,\text{ and }
\|\nabla^2 q_j\|_2\le1,\quad\text{for all }j\in J_2,
\end{equation}
which can always be achieved by appropriate renorming of $q_j$. (Here $\nabla^2 q_j$ denotes the 
(constant) Hesse matrix of $q_j$.)

\begin{lemma} \label{Hardyl}
There is a constant $K$ depending only on $\Omega$, the choice of $\Omega_j$, $j\in J_2$, and $m\ge1$,
 such that for all $j\in J_2$ and $u\in H^m(\Omega)\cap H^1_0(\Omega)$, 
\begin{equation} \label{Hardye}
|u/q_j|_{H^{m-1}(\Omega_j)}\le K \|u\|_{H^{m}(\Omega_j)}.
\end{equation}  
\end{lemma} 
\pf 
The lemma can be proved following the lines of the proof of 
\cite[Theorem~6.1]{HRW01}. Indeed, $q_j$ is a smooth function and the only difference to the setting of
that theorem is that $\Gamma_j$ is a proper part of the boundary of $\Omega_j$ rather than the full
boundary. The argumentation in the proof however remains valid for this case. 
\eop

We define an operator %
$I_\tri: H^3(\Omega)\cap H^1_0(\Omega)\to S_{d,0}(\tri)$ of interpolatory type. 
For all $T\in\tri_0\cup\tri_P$ we set $I_\tri u|_T=I_T(u|_T)$, with the local operators 
$I_T$ defined as follows.

If $T\in \tri_0$, then $I_Tu$ is the polynomial of degree $d$ that
interpolates $u$ on $D_{d,T}$, that is $I_Tu\in\mathbb{P}_d$ is uniquely determined by the conditions
$I_Tu(\xi)=u(\xi)$ for all $\xi\in D_{d,T}$, see e.g.~\cite[Theorem 1.11]{LSbook}.

Let $T\in\tri_P$ with the curved edge on $\Gamma_j$. Denote by $T^{**}$ the (unique) triangle with 
straight edges such that $B_T\subset T^{**}\subset T$. (Note that $T^{**}=T^*$ if $T$ is convex.)
Then  $I_Tu:=pq_j$, where $p\in\mathbb{P}_{d-1}$ interpolates $u/q_j$ on $D_{d-1,T^{**}}$. This
interpolation scheme is well defined even though $T^{**}$ may include points on the boundary
$\Gamma_j$ because $u/q_j\in H^2(\Omega_j)$ by Lemma~\ref{Hardyl}, and hence $u/q_j$ can be identified with a
continuous function on $\overline{\Omega}_j$ by Sobolev embedding.

Finally, assume that $T\in\tri_B$. Then $p:=I_\tri u|_T\in\mathbb{P}_{d+1}$
is determined by the following interpolation conditions on $D_{d+1,T}$: $p(\xi)=u(\xi)$ if 
$\xi\in D^0_{d+1,T}$, and $p(\xi)=I_{T'}u(\xi)$ if $\xi\in D_{d+1,T}\setminus D^0_{d+1,T}$, where 
$T'$ is a triangle in $\tri_0\cup\tri_P$ containing the point $\xi$. The triangle $T'$ is uniquely
defined unless $\xi$ is the interior vertex of a pie-shaped triangle. In the latter case it is easy
to check that $p(\xi)=u(\xi)$ independently of the choice of $T'$, which shows that $p$ is 
well defined. This argument also shows that $I_\tri u$ is continuous at the interior vertices of all
pie-shaped triangles.

To see that $I_\tri u\in S_{d,0}(\tri)$ we need to demonstrate the continuity of $I_\tri u$ across all interior
edges of $\tri$. If $e$ is the common edge of two triangles $T',T''\in\tri_0$, then the continuity follows from
the standard argument that the restrictions $I_{T'}u|_e$ and  $I_{T''}u|_e$ coincide as univariate 
polynomials of degree $d$ satisfying identical interpolation conditions at $d+1$ points of
$D_{d,T'}\cap e=D_{d,T''}\cap e$. The same argument applies to the common edges of buffer triangles.
Finally, the continuity of $I_\tri u$ across edges shared by buffer triangles with either ordinary or
pie-shaped triangles follows from the exact reproduction of univariate polynomials of degree at most 
$d+1$ by the interpolation polynomial $I_\tri u|_T$ on the edges of a buffer triangle $T$.

\begin{lemma} \label{pie}
Let $T\in\tri_P$ and its curved edge $e\subset\Gamma_j$. Then
\begin{equation}\label{pieb1}
\|I_Tu\|_{L^\infty(T)}\le C_1h_T\|u/q_j\|_{L^\infty(T)}\quad\text{if }\;u\in H^3(\Omega)\cap H^1_0(\Omega),
\end{equation}
where $C_1$ depends only on   $d$ and $h_T/\rho_T$.
Moreover, if $d\ge2$ and $3\le m\le d+1$, then for any $u\in H^m(\Omega)\cap H^1_0(\Omega)$, 
\begin{align}\label{pieb2}
\|u-I_{T}u\|_{H^k(T)}&\le C_2h_T^{m-k}|u/q_j|_{H^{m-1}(T)},\quad k=0,\ldots,m-1,\\
\label{pieb3}
\|u-I_{T}u\|_{L^\infty(T)}&\le C_3h_T^{m-1}|u/q_j|_{H^{m-1}(T)},
\end{align}
where $C_2,C_3$ depend only on   $d$ and $h_T/\rho_T$.
\end{lemma} 

\pf We will denote by $C$ various constants depending only on $d$, and by $\tilde C$
constants depending only on $d$ and $h_T/\rho_T$.

Assume that $u\in H^3(\Omega)\cap H^1_0(\Omega)$ and recall that by definition $I_Tu=pq_j$, 
where $p\in\mathbb{P}_{d-1}$ interpolates $u/q_j$ on $D_{d-1,T^{**}}$.
It is not difficult to derive from the proof of \cite[Theorem 1.12]{LSbook} that
\begin{equation}\label{intb}
\|p\|_{L^\infty(T^{**})}\le C\max_{\xi\in D_{d-1,T^{**}}}|p(\xi)|,
\end{equation}
which implies
$$
\|p\|_{L^\infty(T^{**})}\le C \|u/q_j\|_{L^\infty(T^{**})}\le C \|u/q_j\|_{L^\infty(T)}.$$
Since $T$ is star-shaped with respect to $B_T\subset T^{**}$, for any $x\in T$ the absolute value of the 
restriction of $p$ to the straight
line through $x$ and the incenter $z_T$ of $T^{**}$ is bounded by $\|p\|_{L^\infty(T^{**})}$ inside 
$B_T$. By the well-known extremal property of the Chebyshev polynomial $T_{d-1}$ it follows that 
$|p(x)|\le |T_{d-1}(\|x-z_T\|_2/\rho_T)|\,\|p\|_{L^\infty(T^{**})}$, so that
$$
\|p\|_{L^\infty(T)}\le\tilde C\|p\|_{L^\infty(T^{**})},$$
with $\tilde C=|T_{d-1}(h_T/\rho_T)|$. In view of \eqref{q0}, $\|q_j\|_{L^\infty(T)}\le h_T$, and hence
$$
\|I_Tu\|_{L^\infty(T)}=\|pq_j\|_{L^\infty(T)}\le h_T\|p\|_{L^\infty(T)},$$
which completes the proof of \eqref{pieb1}.

Since the area of $T$ is less or equal $\frac{\pi}{4}h^2_T$ and $\p^\alpha (I_Tu)\in\mathbb{P}_{d-k+1}$
if $|\alpha|=k$, it follows that
$$
\|\p^\alpha (I_Tu)\|_{L^2(T)}\le \frac{\sqrt{\pi}}{2}h_T\|\p^\alpha (I_Tu)\|_{L^\infty(T)}
\le \tilde Ch_T\|\p^\alpha (I_Tu)\|_{L^\infty(T^{**})}.$$
By Markov inequality (see e.g.~\cite[Theorem 1.2]{LSbook}) we get furthermore
$$
\|\p^\alpha (I_Tu)\|_{L^\infty(T^{**})}\le
\frac{C}{\rho_T^k}\|I_Tu\|_{L^\infty(T^{**})},$$
and hence
\begin{equation}\label{Markov}
|I_Tu|_{H^k(T)}\le \tilde C h^{1-k}_T\|I_Tu\|_{L^\infty(T^{**})}.
\end{equation}
In view of  \eqref{pieb1} we  arrive at
\begin{equation}\label{pieb4}
|I_Tu|_{H^k(T)}\le \tilde Ch^{2-k}_T\|u/q_j\|_{L^\infty(T)},\quad \text{if }\;u\in
H^3(\Omega)\cap H^1_0(\Omega).
\end{equation}

Let $d\ge2$ and $3\le m\le d+1$, and let $u\in H^m(\Omega)\cap H^1_0(\Omega)$. It follows from
Lemma~\ref{Hardyl} that $u/q_j\in H^{m-1}(T)$.
By the results in \cite[Chapter 4]{BrennerScott} there exists a polynomial $\tilde p\in\mathbb{P}_{m-2}$ such that 
\begin{align*}
\|u/q_j-\tilde p\|_{H^k(T)}&\le \tilde Ch_T^{m-k-1}|u/q_j|_{H^{m-1}(T)},\quad k=0,\ldots,m-1,\quad\text{and}\\
\|u/q_j-\tilde p\|_{L^\infty(T)}&\le \tilde Ch_T^{m-2}|u/q_j|_{H^{m-1}(T)}.
\end{align*}
Indeed, a suitable $\tilde p$ is given by the
\emph{averaged Taylor polynomial} \cite[Definition 4.1.3]{BrennerScott} with respect to the disk $B_T$,
and the inequalities in the last display follow from \cite[Lemma 4.3.8]{BrennerScott} (Bramble-Hilbert Lemma) and
\cite[Proposition 4.3.2]{BrennerScott}, respectively. (It is easy to check by inspecting the proofs in 
\cite{BrennerScott} that the quotient $d_T/\rho_T$ can be used in the estimates instead of the so-called
chunkiness parameter used there.)

Since
$$
u-I_Tu=(u/q_j-\tilde p)q_j-I_T(u - \tilde pq_j), $$
we have for any norm $\|\cdot\|$,
$$
\|u-I_Tu\|\le\|(u/q_j-\tilde p)q_j\|+\|I_T(u - \tilde pq_j)\|. $$
In view of \eqref{q0}, 
$$
\|(u/q_j-\tilde p)q_j\|_{L^\infty(T)}\le h_T \|u/q_j-\tilde p\|_{L^\infty(T)}
\le \tilde Ch_T^{m-1}|u/q_j|_{H^{m-1}(T)},$$
and for any $k=0,\ldots,m-1$,
\begin{align*}
\|(u/q_j-\tilde p)q_j\|_{H^k(T)}
&\le Ch_T\|u/q_j-\tilde p\|_{H^k(T)}+C\|u/q_j-\tilde p\|_{H^{k-1}(T)})\\
&\le \tilde Ch_T^{m-k}|u/q_j|_{H^{m-1}(T)}.
\end{align*}
Furthermore, by \eqref{pieb1}
$$
\|I_T(u - \tilde pq_j)\|_{L^\infty(T)}\le \tilde C h_T \|u/q_j - \tilde p\|_{L^\infty(T)}
\le \tilde Ch_T^{m-1}|u/q_j|_{H^{m-1}(T)},$$
and by \eqref{pieb4}
$$
\|I_T(u - \tilde pq_j)\|_{H^k(T)}\le \tilde Ch^{2-k}_T\|u/q_j - \tilde p\|_{L^\infty(T)}
\le \tilde Ch^{m-k}_T|u/q_j|_{H^{m-1}(T)}.$$
By combining the inequalities in the five last displays we deduce \eqref{pieb2} and \eqref{pieb3}.
\eop

The approximation power of the space $S_{d,0}(\tri)$ is estimated in the following theorem.
\begin{theorem} \label{approx}\mpr{aprox}
Let $d\ge2$ and $3\le m\le d+1$. For any $u\in H^m(\Omega)\cap H^1_0(\Omega)$,
\begin{equation}\label{appr1}
\inf_{s\in S_{d,0}(\tri)}\|u-s\|_{H^k(\Omega)}\le
   C_4h^{m-k}\|u\|_{H^m(\Omega)},\quad k=0,\ldots,m-1,
\end{equation}
where $h$ is the maximum diameter of the triangles in $\tri$,
and $C_4$ is a constant depending only on $d$ and $R$, as well as on $\Omega$ and the choice of $\Omega_j$ 
for all $j\in J_2$.
\end{theorem}

\pf We again use $C$ for various constants depending only on the parameters mentioned for $C_4$ 
in the formulation of the theorem.
We show the error bound in \eqref{appr1} for $s=I_\tri u$. For this sake we estimate
the norms of $u-I_Tu$ on all triangles $T\in\tri$. 

If $T\in\tri_0$, then by \cite[Theorem 4.4.4]{BrennerScott} 
\begin{equation}\label{bound0}
\|u-I_{T}u\|_{H^k(T)}\le Ch_T^{m-k}|u|_{H^m(T)},\quad k=0,\ldots,m,
\end{equation}
where $C$ depends only on $d$ and $h_T/\rho_T$.
If $T\in\tri_P$, with the curved edge $e\subset\Gamma_j$, %
then the estimate \eqref{pieb2}  holds by Lemma~\ref{pie}.

Let $T\in\tri_B$, $p:=I_\tri u|_T$ and let $\tilde p\in\mathbb{P}_{d+1}$ be the interpolation polynomial that satisfies
$\tilde p(\xi)=u(\xi)$ for all $\xi\in D_{d+1,T}$. Then 
$$
\tilde p(\xi)-p(\xi)=
\begin{cases} 0 & \text{if }\xi\in D^0_{d+1,T},\\
u(\xi)-I_{T'}u(\xi) & \text{if }\xi\in D_{d+1,T}\setminus D^0_{d+1,T},
\end{cases}$$
where  $T'\in \tri_0\cup\tri_P$  contains $\xi$. Hence, by the same arguments leading to \eqref{intb} and
\eqref{Markov} in the proof of Lemma~\ref{pie}, we conclude that for $k=0,\ldots,m$,
\begin{align*}
\|\tilde p-p\|_{H^k(T)}&\le Ch_T^{1-k}\|\tilde p-p\|_{L^\infty(T)}\\
&\le
 Ch_T^{1-k}\max\{\|u-I_{T'}u\|_{L^\infty(T')}:T'\in \tri_0\cup\tri_P,\; T'\cap T\ne\emptyset\},
\end{align*}
whereas by \cite[Theorem 4.4.4]{BrennerScott} we have
$$
\|u-\tilde p\|_{H^k(T)}\le Ch_T^{m-k}|u|_{H^m(T)},$$
with the constants depending only on $d$ and $h_T/\rho_T$.
If $T'\in\tri_0\cup\tri_P$, then by \cite[Corollary 4.4.7]{BrennerScott} and  \eqref{pieb3},
$$
\|u-I_{T'}u\|_{L^\infty(T')}\le Ch_{T'}^{m-1}
\begin{cases} 
|u|_{H^m(T')}&\text{if }T'\in\tri_0,\\
|u/q_j|_{H^{m-1}(T')}&\text{if }T'\in\tri_P,
\end{cases}$$
where $C$ depends only on $d$ and $h_{T'}/\rho_{T'}$. By combining these inequalities we obtain
an estimate of $\|u- p\|_{H^k(T)}$ by $C\tilde h^{m-k}$ times the maximum of $|u|_{H^m(T)}$, 
$|u|_{H^m(T')}$ for $T'\in\tri_0$ sharing edges with $T$, and $|u/q_j|_{H^{m-1}(T')}$
for $T'\in\tri_P$ sharing edges with $T$. Here $C$ depends only on $d$ and the maximum of 
$h_{T}/\rho_{T}$ and $h_{T'}/\rho_{T'}$, and $\tilde h$ is the maximum of $h_T$ and $h_{T'}$ for 
$T'\in\tri_0\cup\tri_P$  sharing edges with $T$.

By using \eqref{bound0} on $T\in\tri_0$, \eqref{pieb2} on $T\in\tri_P$ and the estimate of the last
paragraph on $T\in\tri_B$, we get
\begin{align*}
\|u-I_\tri u\|_{H^k(\Omega)}^2&=\sum_{T\in\tri}\|u-I_\tri u\|_{H^k(T)}^2\\
&\le Ch^{2(m-k)}\Big(\sum_{T\in\tri_0\cup\tri_B}|u|_{H^m(T)}^2
+\sum_{T\in\tri_P}|u/q_{j(T)}|_{H^{m-1}(T)}^2\Big),
\end{align*}
where $C$ depends only on $d$ and $R$, and $j(T)$ is the index of $\Gamma_j$ containing the curved edge
of $T\in\tri_P$. Clearly,
$$
\sum_{T\in\tri_0\cup\tri_B}|u|_{H^m(T)}^2\le|u|_{H^m(\Omega)}^2,$$
whereas by Lemma~\ref{Hardyl},
$$
\sum_{T\in\tri_P}|u/q_{j(T)}|_{H^{m-1}(T)}^2\le\sum_{j\in J_2}|u/q_{j}|_{H^{m-1}(\Omega_j)}^2
\le K \|u\|_{H^m(\Omega)}^2,$$
where $K$ is the constant of \eqref{Hardye} depending only on $\Omega$ and the choice of $\Omega_j$, 
 $j\in J_2$.
\eop

\subsection{Solution of elliptic problems}
For simplicity we restrict ourselves to symmetric regular elliptic problems with 
homogeneous Dirichlet boundary conditions.

Consider the  variational problem
\begin{equation}
\label{coercive}
\text{find $u\in H^1_0(\Omega)$ such that 
$a(u,v)=(f,v)$ for all $v\in H^1_0(\Omega)$,}
\end{equation}
where $a(\cdot,\cdot)$ is a \emph{symmetric}, \emph{bounded} and \emph{coercive} 
bilinear form (see e.g.~\cite{BrennerScott} for the definitions), 
such that the \emph{regularity condition}
\begin{equation}
\label{regular}
\|u\|_{H^2(\Omega)}\le C_R\|f\|_{L^2(\Omega)},
\end{equation}
holds with some constant $C_R>0$ independent of $f$.
Note that  (\ref{regular}) holds for the domains
considered in this work if $\omega_j\le\pi$, $j=1,\ldots,n$, see e.g.~\cite[p.~158]{Schwab98}. 

The finite element approximation $\tilde u$  of $u$ relying on the space $S_{d,0}(\tri)$ 
is found as the solution of the discretized problem
\begin{equation}
\label{discr1}
\text{find $\uw\in S_{d,0}(\tri)$ such that 
$a(\uw,s)=(f,s)$ for all $s\in S_{d,0}(\tri)$.}
\end{equation}
The existence and uniqueness of the solutions of (\ref{coercive}) and \eqref{discr1}  follows by the
well-known Lax-Milgram Theorem as soon as the bilinear form is coercive and bounded.

As a consequence of Theorem~\ref{approx}, we obtain the following 
error estimates by using the standard arguments in the finite element method, see the proofs of 
Theorems~5.4.4, 5.4.8 and 5.7.6 in \cite{BrennerScott}.

\begin{theorem}\label{bound1}\mpr{bound1}
Suppose that the variational problem (\ref{coercive}) is symmetric, bounded, coercive and regular
in the sense of (\ref{regular}). Let $d\ge2$ and $3\le m\le d+1$. Then 
\begin{eqnarray}
&&\|u-\uw\|_{L^2(\Omega)}\le
   C_1h^m\|u\|_{H^m(\Omega)},\label{estDir1}\\
&&\|u-\uw\|_{H^1(\Omega)}\le
   C_2h^{m-1}\|u\|_{H^m(\Omega)},\label{estDir2}
\end{eqnarray}
where $h$ is the maximum diameter of the triangles in $\tri$,
and $C_1,C_2$ are constants depending only on $\Omega$, $d$, $R$ and
the choice of $\Omega_j$, $j\in J_2$.
\end{theorem}

\section{Bernstein-B\'ezier basis}

In order to treat elliptic problems with homogenous Dirichlet boundary conditions,
we need to describe a suitable local basis for $S_{d,0}(\tri)$.
To this end we use the Bernstein-B\'ezier techniques and construct a stable 
minimal determining set (MDS) for  $S_{d,0}(\tri)$.
The main idea is to factorize polynomials over pie-shaped triangles.
As usual, we identify the functionals of the MDS with certain domain points, 
albeit not always in the standard way.

\subsection{Bernstein polynomials}
Let $T$ be a non-degenerate triangle in the plane with vertices $v_{1},v_{2},v_{3}$. The 
bivariate Bernstein polynomials with respect to $T$ are defined by 
 \begin{equation*}
B_{ijk}^{d}(v):= \frac{d!}{i!j!k!}b_{1}^{i}b_{2}^{j}b_{3}^{k},\quad i+j+k=d,  
\end{equation*}
where $b_{1},b_{2},b_{3}$ are the barycentric coordinates of $v$, that is the unique coefficients of the
expansion $v=\sum _{i=1}^{3}b_{i}v_{i}$ with $\sum _{i=1}^{3}b_{i}=1$. Recall that Bernstein polynomials
form a basis for $\mathbb{P}_d$ and have many useful properties for dealing with bivariate polynomials
and piecewise polynomials, see \cite{LSbook}. In particular, every polynomial $p$ of degree $d$ 
can be written in the \emph{BB-form} as
\begin{equation}\label{BBform}
p=\sum _{i+j+k=d}c_{ijk}B_{ijk}^{d},
\end{equation}
where $c_{ijk}$ are called the \emph{BB-coefficients} of $p$.

The BB-form is stable in the sense that there is a constant $K>0$ depending
only on $d$, such that
\begin{equation}\label{BBstab}
K\|c\|_{\infty}\le\|p\|_{L^\infty(T)}\le \|c\|_{\infty},
\qquad \|c\|_{\infty}:=\max_{i+j+k=d}|c_{ijk}|.
\end{equation}

It is often convenient to index Bernstein polynomials by the elements of the set
\begin{equation}\label{domainpoints}
D_{d,T}:=\left\lbrace \xi _{ijk} =\frac{iv_{1}+jv_{2}+kv_{3}}{d}:\;i+j+k=d,\;i,j,k \geq 0 \right\rbrace
\end{equation}
of so-called  \emph{domain points}, such that $B_{\xi}^{d}:=B_{ijk}^{d}$ and $c_{\xi}:=c_{ijk}$ when
$\xi=\xi_{ijk}\in D_{d,T}$.

In particular, it is easy to express the continuity of piecewise polynomials as follows. 
Given two triangles $T$ and $\tilde T$ sharing an edge $e$, let $p$ and $\tilde p$
be two polynomials of degree $d$ written in the BB-form
$$
p=\sum _{\xi\in D_{d,T}}c_{\xi}B_{\xi}^{d}\quad\text{ and }\quad
\tilde p=\sum _{\xi\in D_{d,\tilde T}}\tilde c_{\xi}\tilde B_{\xi}^{d},$$
where $B_{\xi}^{d}$ and $\tilde B_{\xi}^{d}$ are the Bernstein polynomials with respect to $T$ and
$\tilde T$, respectively. Then $p$ and $\tilde p$ join  continuously along $e$ if their 
BB-coefficients over $e$ coincide, i.e. 
\begin{equation}\label{smoothnessc0}
\tilde{c}_{\xi}=c_{\xi},\quad\text{ for all }\; \xi\in D_{d,T}\cap D_{d,\tilde T}
\end{equation}

\subsection{Minimal determining sets}

A key concept for dealing with spline spaces using 
Bernstein-B\'ezier techniques
is that of a \emph{stable local minimal determining set of domain points}, see \cite{LSbook}. 
Our setting of curved domains requires a bit more general concept of a minimal
determining set of \emph{functionals} which we describe below.

\begin{definition}\label{MDS}
\rm
Let $S$ be a finite dimensional linear space and $S^*$ its dual space. 
A set $\Lambda \subset S^*$ is said to be a \emph{determining set} for $S$ if 
for any $s\in S$
\begin{equation*}\label{e7}
\lambda(s)=0 \quad\forall  \lambda \in \Lambda\quad 
\Longrightarrow \quad s=0,
\end{equation*}
and $\Lambda$ is a \emph{minimal determining set (MDS)} for the space $S$ if there is no smaller determining
set. 
\end{definition}

In other words, a determining set is a spanning set of linear functionals, and
an MDS is a basis of $S^*$. Therefore, any MDS $\Lambda$ uniquely 
determines a \emph{basis} 
\begin{equation}\label{basis}
\{s_\lambda:\lambda\in\Lambda\}
\end{equation}
 of $S$ by duality, 
such that
$$
\lambda(s_\mu)=\delta_{\lambda,\mu},\qquad \lambda,\mu\in\Lambda,$$
and any spline $s\in S$ can be written uniquely in the form
$$
s=\sum_{\lambda\in\Lambda}c_\lambda s_\lambda,\quad\text{with }
c_\lambda=\lambda(s)\in\mathbb{R}.$$
A minimal determining set serves as a set of \emph{degrees of freedom} of the
space $S$, as the value $\mu(s)$ of any functional $\mu\in S^*$ 
is uniquely determined by the values $\lambda(s)$ for all
$\lambda\in\Lambda$. For example, if $S$ is a space of functions on a set
$\Omega$ such that the point evaluation functionals $s(x)$, $x\in\Omega$,
are well-defined, then the values $s(x)$ for all $x$ are determined by 
$\lambda(s)$, $\lambda\in\Lambda$. 

We define the properties of \emph{locality} and \emph{stability} of an MDS 
for the spaces of piecewise polynomial splines. 

Let $\tri$ be a \emph{partition} of a bounded domain $\Omega\subset\mathbb{R}^n$ 
into a finite number of \emph{cells}
$T\in\tri$, where each $T$ is a closed set with dense interior, 
such that $\bigcup_{T\in\tri}T=\overline{\Omega}$ and 
$T_1\cap T_2\subset\p T_1\cap \p T_2$ for any $T_1, T_2\in\tri$.
A finite dimensional linear space $S$ of real valued functions defined on $\Omega$ is said to be a
\emph{spline space} with respect to a partition $\tri$ if 
$s|_T$ is an $n$-variate algebraic polynomial for any $s\in S$ and any $T\in\tri$. 
The maximal total degree of such polynomials is called the \emph{degree} of the spline space and
is denoted $d_S$.

The \emph{star} of a set $A\subset\Omega$ with respect to a partition $\tri$
is given by
$$
\st(A):=\bigcup\{T\in\tri: T\cap A\ne\emptyset\}.$$
For any positive integer $\ell$ we define the \emph{$\ell$-star} of $A$ recursively as
$$
\st^1(A):=\st(A),\quad\st^\ell(A):=\st(\st^{\ell-1}(A)),\quad \ell\ge2.$$
It is easy to see that $B\cap\st^\ell(A)\ne\emptyset$ if and only if there
exists a chain of cells $T_1,\ldots,T_\ell\in\tri$ such that 
$T_1\cap A\ne\emptyset$, $T_\ell\cap B\ne\emptyset$, and 
$T_i\cap T_{i+1}\ne\emptyset$ for all $i=1,\ldots,\ell-1$.
It is easy to derive from this that
\begin{align*}
B\cap\st^\ell(A)\ne\emptyset\quad&\Longleftrightarrow\quad A\cap\st^\ell(B)\ne\emptyset,
\qquad A,B \subset\Omega,\\
T\subset\st^\ell(A)\quad&\Longleftrightarrow\quad A\cap\st^{\ell-1}(T)\ne\emptyset,
\qquad T\in\tri,\quad A \subset\Omega,\\
T_2\subset\st^\ell(T_1)\quad&\Longleftrightarrow\quad T_1\subset\st^\ell(T_2),
\qquad T_1,T_2\in\tri.
\end{align*}

Given a spline space $S$, a set $\omega\subset\Omega$ is said to be a
\emph{supporting set} of a linear functional $\lambda\in S^*$ if 
\begin{equation}\label{supset}
\lambda(s)=0\quad\text{for all $s$ such that $s|_\omega=0$}.
\end{equation}
Note that a supporting set is not unique, and a minimal supporting set may not
exists since the intersection of two supporting sets is not necessarily a
supporting set. For example, if $S$ is a space of continuously differentiable splines on
a triangulation $\tri$, and $\lambda(s)$ is a partial derivative of $s$ evaluated at a
vertex $v$ of $\tri$, then every triangle attached to $v$ is a supporting set of
$\lambda$, but their intersection is just $\{v\}$ and this is in general not a
supporting set.

Given a spline space $S$ and an MDS $\Lambda$, we define for each
$T\in\tri$ the set
$$
\Lambda_T:=\{\lambda\in\Lambda: T\subset\supp s_\lambda\},$$
where $\{s_\lambda:\lambda\in\Lambda\}$ is the basis of $S$ dual to $\Lambda$.
Thus, $\lambda\in\Lambda_T$ if and only if for a spline $s\in S$, 
$s|_T$ depends on the coefficient $c_\lambda=\lambda(s)$. The number
$$
\kappa_\Lambda:=\max_{T\in\tri}|\Lambda_T|$$
is called the \emph{covering number} of the MDS $\Lambda$.

\begin{definition}\label{localMDS}
\rm
A minimal determining set $\Lambda$ for a spline space $S$ is said to be
\emph{$\ell$-local} if there is a family of supporting sets $\omega_\lambda$ of 
$\lambda\in\Lambda$ such that $\omega_\lambda\subset \st^{\ell}(T)$
for any $T\in\tri$ such that $\lambda\in\Lambda_T$. 
\end{definition}

\begin{lemma}\label{locbasis}
Let $\Lambda$ be an $\ell$-local MDS. Then for any $T\subset\supp s_\lambda$,
$$
\omega_\lambda\subset \st^{\ell}(T)\quad\text{\rm and}\quad
\supp s_\lambda\subset\st^{2\ell+1}(T).$$
Moreover, if $\omega_\lambda=T_\lambda\in\tri$, then
$$
\supp s_\lambda\subset\st^{\ell}(T_\lambda).$$
\end{lemma}
\pf 
Clearly, the support of any spline, in particular $s_\lambda$, is a union of
cells $T\in\tri$. Let $T\subset\supp s_\lambda$. Then $\lambda\in \Lambda_T$ and
hence $\omega_\lambda\subset \st^{\ell}(T)$. 
Therefore $T\subset\st^{\ell+1}(\omega_\lambda)$. It follows that
$\supp s_\lambda\subset\st^{\ell+1}(\omega_\lambda)\subset\st^{2\ell+1}(T)$.
If $\omega_\lambda=T_\lambda\in\tri$, then $T_\lambda\subset \st^{\ell}(T)$
implies $T\subset\st^{\ell}(T_\lambda)$, and hence 
$\supp s_\lambda\subset\st^{\ell}(T_\lambda)$.
 \eop
 
As a consequence of this lemma we note that $\kappa_\Lambda$ is bounded by the product
of the dimension of the polynomials of degree $d_S$ and the maximal possible number of
cells in $\st^{2\ell+1}(T)$. Indeed, this product is an upper bound for the dimension
of $S|_{\st^{2\ell+1}(T)}$, and hence for $|\Lambda_T|$ because all basis splines
$s_\lambda$ with $\lambda\in\Lambda_T$ are zero outside of $\st^{2\ell+1}(T)$ and
therefore must be linearly independent within this set.
 
We now assume that $\Lambda$ is an $\ell$-local MDS for a spline space $S$, 
and $\{\omega_\lambda:\lambda\in \Lambda\}$ is the corresponding family of 
supporting sets. According to \eqref{supset}, $\lambda(s)$ is completely
determined by the restricted spline $s|_{\omega_\lambda}$, which means that it can be
considered as a bounded linear functional on the space 
$S|_{\omega_\lambda}=\{s|_{\omega_\lambda}:s\in S\}$, and so there is a
constant $K_1\in(0,\infty)$ such that
\begin{equation}\label{K1}
|\lambda(s)|\le K_1\|s\|_{L^\infty(\omega_\lambda)}\quad
\text{for all }\lambda\in\Lambda,\; s\in S.
\end{equation}
In view of the definition of $\Lambda_T$, $s|_T=0$ if $\lambda(s)=0$ for all
$\lambda\in\Lambda_T$. Hence $\max_{\lambda\in\Lambda_T}|\lambda(s)|$ is a norm
on the space $S|_T$ and, in view of the equivalence of all norms on a
finite-dimensional space, there exists a constant $K_2\in(0,\infty)$ such that
\begin{equation}\label{K2}
\|s\|_{L^\infty(T)}\le K_2\max_{\lambda\in\Lambda_T}|\lambda(s)|\quad
\text{for all }T\in\tri,\; s\in S.
\end{equation}
We call any constants $K_1,K_2$ satisfying \eqref{K1} and \eqref{K2}, respectively, 
the \emph{stability constants} of the MDS $\Lambda$.

\begin{proposition}\label{Lpstab}
Let $\Lambda$ be an MDS for a spline space $S$. Then for any 
$a=\{a_\lambda\}_{\lambda\in\Lambda}\subset\mathbb{R}$,
$$
K_1^{-1}\|a\|_\infty\le\|\sum_{\lambda\in\Lambda}a_\lambda s_\lambda\|_{L^\infty(\Omega)}
\le \kappa_\Lambda K_2\|a\|_\infty,$$
with $\|a\|_\infty:=\max_{\lambda\in\Lambda}|a_\lambda|$.
Moreover, for any $1\le p<\infty$,
$$
\|\sum_{\lambda\in\Lambda}a_\lambda s^{(p)}_\lambda\|_{L^p(\Omega)}
\le \kappa_\Lambda^{1-1/p} K_2\|a\|_p$$
where 
$$
s^{(p)}_\lambda:=|\supp s_\lambda|^{-1/p}s_\lambda,\qquad
\|a\|_p:=\Big(\sum_{\lambda\in\Lambda}|a_\lambda|^p\Big)^{1/p}.$$
\end{proposition}
\pf 
Let $s=\sum_{\lambda\in\Lambda}a_\lambda s_\lambda$. By \eqref{K1},
$|a_\lambda|=|\lambda(s)|\le K_1\|s\|_{L^\infty(\Omega)}$, which implies the first inequality.
On the other hand,  for any $T\in\tri$, 
since $\|s_\lambda\|_{L^\infty(\Omega)}\le K_2$ by \eqref{K2},
$$
\|s|_T\|_{L^\infty(T)}=\|\sum_{\lambda\in\Lambda_T}a_\lambda s_\lambda\|_{L^\infty(T)}
\le \kappa_\Lambda\max_{\lambda\in\Lambda_T}
\|a_\lambda s_\lambda\|_{L^\infty(T)}\le\kappa_\Lambda K_2\|a\|_\infty,$$
which completes the proof in the case $p=\infty$.

Let now $1\le p<\infty$ and $s=\sum_{\lambda\in\Lambda}a_\lambda s^{(p)}_\lambda$. 
Then for any $T\in\tri$,
$$
\|s|_T\|^p_{L^p(T)}=\|\sum_{\lambda\in\Lambda_T}a_\lambda s^{(p)}_\lambda\|^p_{L^p(T)}
\le \kappa_\Lambda^{p-1}\sum_{\lambda\in\Lambda_T}
\|a_\lambda s^{(p)}_\lambda\|^p_{L^p(T)}.$$ 
Hence
$$
\|s\|^p_{L^p(\Omega)}=\sum_{T\in\tri}\|s|_T\|^p_{L^p(T)}\le
\kappa_\Lambda^{p-1}\sum_{\lambda\in\Lambda}|a_\lambda|^p
\sum_{T\in \tri\atop T\subset\supp s_\lambda}\|s^{(p)}_\lambda\|^p_{L^p(T)}
\le\kappa_\Lambda^{p-1}K_2\|a\|^p_p,$$
where we have taken into account that 
$\|s_\lambda\|_{L^p(\Omega)}\le K_2 |\supp s_\lambda|^{1/p}$ in view of \eqref{K2}.
This implies the second statement.
\eop

This proposition shows that $\kappa_\Lambda K_1 K_2$ is an upper bound for the
$L_\infty$-stability  constant of the basis $\{s_\lambda\}_{\lambda\in\Lambda}$. 
For $p<\infty$, the lower $L_p$-stability estimates of the form
$$
CK_1^{-1}\|a\|_p\le\|\sum_{\lambda\in\Lambda}a_\lambda s^{(p)}_\lambda\|_{L^p(\Omega)}$$
with some $C>0$ are more complicated. Under additional assumptions on the partition $\tri$,
they hold with $C$ depending only on $\ell$, $n$ and $d_S$, see e.g.~\cite[Lemma 6.2]{D01} and
\cite[Theorem 5.22]{LSbook}.

This motivates the following definition.

\begin{definition}\label{stableMDS}
\rm
A family of minimal determining sets $\Lambda_i$, $i\in \mathcal{I}$, 
for given spline spaces $S_i$, $i\in \mathcal{I}$, is said to be \emph{stable and local} 
if there exist $\ell$, $\kappa$, $K_1$ and $K_2$ such that every $\Lambda_i$ is 
$\ell$-local, the covering numbers satisfy $\kappa_{\Lambda_i}\le \kappa$, $i\in \mathcal{I}$, 
and \eqref{K1} and \eqref{K2} hold for all $\Lambda_i$ with the same 
stability constants $K_1,K_2$.
\end{definition}

Typically such a stable and local family of minimal
determining sets is produced by an algorithm that 
generates an MDS for a specific type of splines on arbitrary regular triangulations in 2D or
3D, see \cite{LSbook}. With some abuse of 
terminology we say in this case that
the \emph{MDS} (given by an algorithm) is \emph{stable and local}.
The parameters $\ell,\kappa,K_1,K_2$ usually depend only on the maximum
degree of polynomials in $S|_T$ and on a shape regularity measure of the cells, such as
the minimum angle of the triangles $T\in\tri$ in the 2D case. 

We now consider the special case when $\tri$ is a \emph{triangulation} 
of a bounded polygonal domain  $\Omega \subset \mathbb{R}^2$, that is each cell $T$ is a
triangle, and the intersection of two different triangles is their common vertex 
or edge if not empty. Let 
$$
D_{d,\triangle}:=\bigcup _{T\in \triangle} D_{d,T},$$
where $D_{d,T}$ is the set of domain points  \eqref{domainpoints} for a single triangle.
In view of \eqref{BBform} every \emph{continuous} piecewise polynomial spline $s$ of degree 
$d$ with respect to $\tri$ can be uniquely represented over each $T\in\tri$ as 
\begin{equation*}
s\vert _{T} = \sum _{\xi \in D_{d,T}}c_{\xi}B_{\xi}^{T,d},
\end{equation*}
where $B_{\xi}^{T,d}$ are BB-basis polynomials of degree $d$ associated with the triangle $T$.
In view of \eqref{smoothnessc0},  the
continuity of $s$ implies that the BB-vectors of $s\vert _{T}$ and $s\vert _{\tilde{T}}$ agree on 
domain points on the edge shared by triangles $T$ and $\tilde{T}$. Therefore each domain point
$\xi\in D_{d,\triangle}$ defines on the space $S^{0}_{d}(\triangle )$ of all
continuous splines of degree $d$ a linear functional $\gamma_\xi$ that picks the BB-coefficient $c_\xi$ 
of $s|_T$ for a triangle $T$ containing $\xi$. It is easy to see that
$$
\Gamma:=\{\gamma_\xi:\xi\in D_{d,\triangle}\}$$
is a minimal determining set for $S^{0}_{d}(\triangle)$ \cite[Section 5.4]{LSbook}. Minimal
determining sets for many subspaces of  $S^{0}_{d}(\triangle )$ can be obtained as subsets of 
$\Gamma$, which are conveniently identified with the corresponding subsets of 
$D_{d,\triangle}$, see \cite{LSbook}.

It is easy to see that a supporting set of $\gamma_\xi$ is given by any triangle
$T\in\tri$ that contains $\xi$. However, $\{\xi\}$ is a supporting
set if  $\xi$ is a vertex of $\tri$, and an edge $e$ of $\tri$ is a supporting set if 
$\xi\in e$. Hence, our generalized definition  of a local MDS
(Definition~\ref{localMDS}) reduces to that given in
\cite[Definition~5.16]{LSbook}, and the parameter $\ell$ is the same if we always take 
a triangle $T_\xi$ containing $\xi$ as a supporting set. Furthermore, 
\cite[Theorem 2.6]{LSbook} shows that \eqref{K1} is satisfied for each $\gamma_\xi$ with the
triangle $T_\xi$ as a supporting set, and $K_1$ depending only on $d$, and  \eqref{K2} is
equivalent to the condition \cite[Eq.~(5.13)]{LSbook}. Therefore, our notion of a stable local
MDS (Definition~\ref{stableMDS}) coincides with \cite[Definition~5.16]{LSbook} in  the case
when the linear functionals belong to $\Gamma$. Similarly, \emph{stable local nodal minimal determining
sets} of  \cite[Section 5.9]{LSbook} are special cases of stable local MDS as defined in 
Definition~\ref{stableMDS}.

\subsection{A stable MDS and stable local basis for $S_{d,0}(\tri)$}\label{sMDS}

\begin{figure}
\centering
\begin{tabular}{cc}
\includegraphics[height=0.4\textwidth]{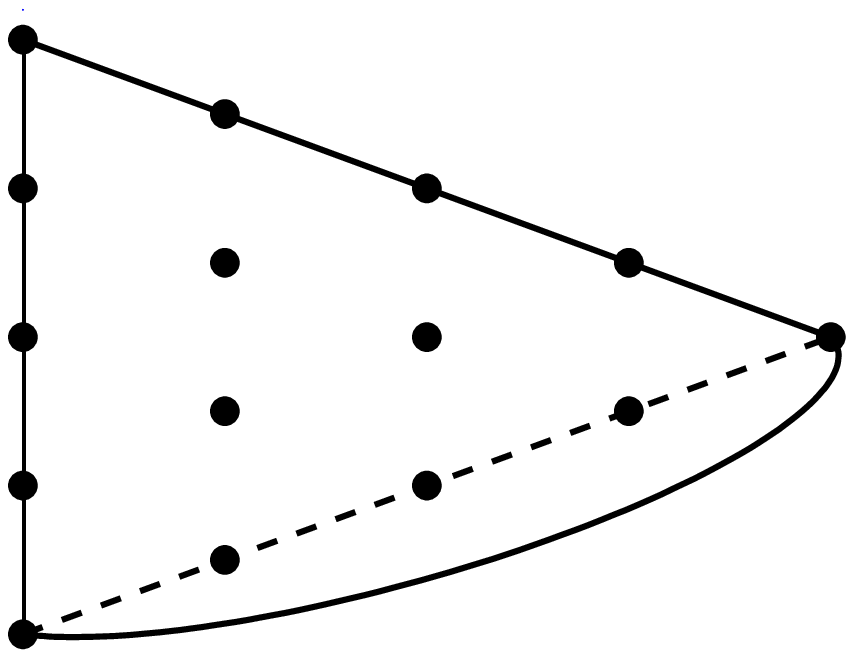} &
\includegraphics[height=0.4\textwidth]{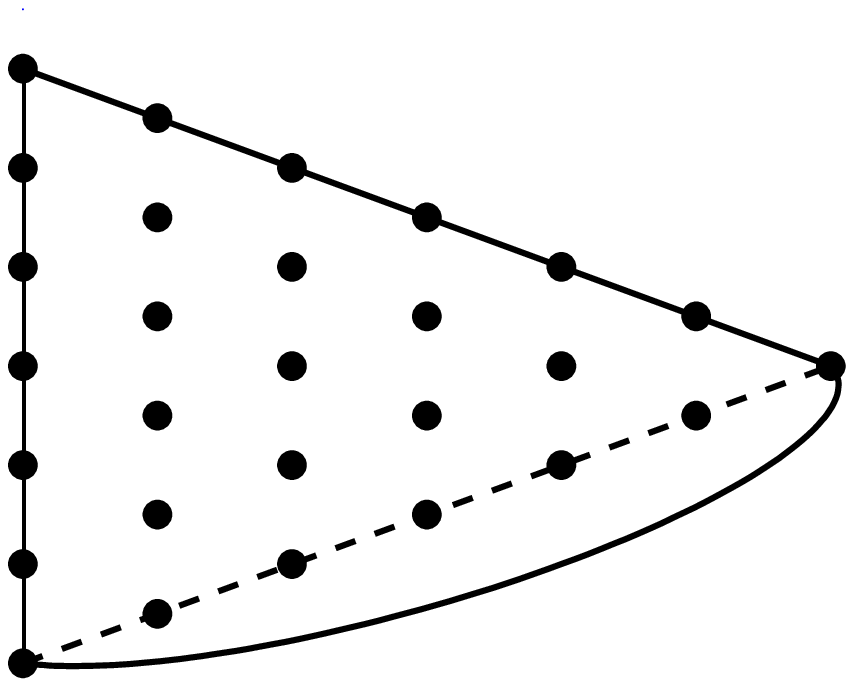} 
\end{tabular}
\caption{
The sets $D^\ast_{d-1,T}$ (left) and $D_{d+1,T^\ast}$ (right)
for $d=5$ over a pie-shaped triangle.
}\label{BezierNetP_C0}
\end{figure}

For each $T \in \tri_P$, with its curved edge $e$ given by the equation $q(x)=0$, 
where the quadratic polynomial is normalized so that
\begin{equation}\label{q5}
q(v)=1\quad\text{for the interior vertex $v$ of $T$,}
\end{equation}
we consider the
space 
$$
q\mathbb{P}_{d-1} =\{pq\sp p\in \mathbb{P}_{d-1}\}
\subset \mathbb{P}_{d+1}.$$
Recall that by B\'ezout theorem every polynomial that vanishes on the conic
$q(x)=0$ is divisible by $q$, and thus $\mathbb{P}_{d-1}q$ consists of
all polynomials of degree $d+1$ that vanish on $e$, that is
$$
q\mathbb{P}_{d-1} =\{ p\in \mathbb{P}_{d+1}\sp p|_{e}=0\}.$$
Since the BB polynomials $B_{ijk}^{d-1}$, $i+j+k=d-1$, w.r.t. $T^\ast$ form a basis for 
$\mathbb{P}_{d-1}$ it is obvious that the set 
$$
\left\lbrace qB_{ijk}^{d-1}\sp i+j+k=d-1\right\rbrace $$ 
is a basis for $q\mathbb{P}_{d-1}$. The set of domain points 
of degree $d-1$ over $T^\ast$ will be denoted $D^\ast_{d-1,T}$. Even though
the set $D^\ast_{d-1,T}$ formally coincides with $D_{d-1,T^\ast}$, the linear 
functionals associated with the domain points are slightly different.
Namely, each $\xi \in D^\ast_{d-1,T}$ represents a linear functional 
$\lambda_\xi$ on $S_{d,0}(\tri)$
which picks the  coefficient $c_\xi$ in the expansion 
$$
s|_T=q\sum_{\xi \in D^\ast_{d-1,T}}c_\xi B_\xi ^{d-1},\qquad 
s\in S_{d,0}(\tri).$$
Clearly, $s|_T\in q\mathbb{P}_{d-1}\subset\mathbb{P}_{d+1}$ may also be expressed 
in terms of BB polynomials $B_\xi ^{d+1}$, $\xi\in D_{d+1,T^\ast}$, 
of degree $d+1$. Figure~\ref{BezierNetP_C0} depicts both sets of domain points
for a pie-shaped triangle $T$ in the case $d=5$.

For the ordinary triangles $T\in\tri_0$ we consider the usual sets of domain
points $D_{d,T}$ (and corresponding functionals $\lambda_\xi=\gamma_\xi$) whose union $D_{d,\tri_0}=\cup_{T\in\tri_0}D_{d,T}$ forms the
standard set of domain points associated with the subtriangulation $\tri_0$ of  
$\tri$.

\begin{figure}
\centering
\begin{tabular}{c}
\includegraphics[height=0.4\textwidth]{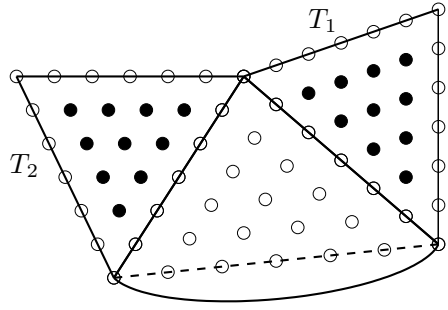} 
\end{tabular}
\caption{
The points in $D^0_{d+1,T_1} \cup D^0_{d+1,T_2}$ for two buffer triangles 
$T_1,T_2\in \tri_{B}$ are marked by black dots ($d=5$).
}\label{BufferMDSC0}
\end{figure}

Finally, for each buffer triangle 
$T:=\left\langle v_1,v_2,v_3\right\rangle\in \tri_{B}$,
let $D_{d+1,T}^0$ be the subset of $D_{d+1,T}$ obtained by removing all domain
points on those edges of $T$ that are shared with triangles in either 
$\tri_0$ or $\tri_P$,
see Figure~\ref{BufferMDSC0}.

We will need the following statement that shows that the division by $q$ on a pie-shaped
triangle is stable.

\begin{lemma}\label{stdiv}
Let $T\in \tri_{P}$ with the curved edge $e$ of $T$
given by the equation $q(x)=0$, where $q=q_j/q_j(v)$ for some $j\in\{1,\ldots,n\}$, and $v$ is
the interior vertex of $T$.
Then for any $r\in\mathbb{P}_{d-1}$,
$$
\|r\|_{L^\infty(T)}\le C\|rq\|_{L^\infty(T)},$$ 
where the constant $C>0$ depends only on $d$ and $B$.
\end{lemma}
\pf
For any
$z\in e$, consider the univariate polynomials $\tilde r(t)=r(tv+(1-t)z)$ and 
$\tilde q(t)=q(tv+(1-t)z)$, $t\in[0,1]$. Then $\tilde q(0)=0$ and 
$\tilde q'(0)=\p_{v-z}q(z)\ge B$ by \eqref{q4}. Hence $\tilde q(t)=t(at+b)$, 
with $b\ge B$. Since $\tilde q(1)=1$, we get $a+b=1$ and hence
$$%
at+b\ge\min\{B,1\},\qquad t\in[0,1].
$$ %
We note for the later use in the proof of Lemma~\ref{Hardyl} that 
$\tilde q'(t)=2at+b$ and hence $\p_{v-z}q(v)=\tilde q'(1)=2a+b=a+1$.

Let $\tilde r(t)=\sum_{i=0}^{d-1}a_it^i$.  In view of the equivalence of all norms on the space of
univariate polynomials of a fixed degree, there exist constants $c_1,c_2>0$ depending only on
$d$ such that
$$
c_1\max_{t\in[0,1]}\Big|\sum_{i=0}^{d-1}a_it^{i}\Big|\le \max_{0\le i\le d-1}|a_i|
\le c_2\max_{t\in[0,1]}\Big|\sum_{i=0}^{d-1}a_it^{i+1}\Big|$$
and hence
$$
\|\tilde r\tilde q\|_{L^\infty(0,1)}\ge 
\min\{B,1\}\max_{t\in[0,1]}\Big|\sum_{i=0}^{d-1}a_it^{i+1}\Big|
\ge \frac{c_1}{c_2}\min\{B,1\}\,\|\tilde r\|_{L^\infty(0,1)}.$$ 
The lemma follows because $T=\bigcup_{z\in e}[z,v]$.  
\eop

We are ready to formulate and prove the main result of this section.
\begin{theorem}\label{S0dMDS}
Let
\begin{equation}
M_0 :=\Big(D_{d,\tri_0}\setminus\p \Omega\Big) \cup 
\bigcup _{T\in \tri_{P}}D_{d-1,T}^\ast \cup
\bigcup _{T\in \tri_{B}} D_{d+1,T}^0.  
\end{equation}
Then $M_0$ is a stable local minimal determining set for the space
$S_{d,0}(\tri)$ with $\ell=1$, $\kappa$ depending only on $d$, $K_1$ depending only on $d$ and $B$, 
and $K_2$ depending only on $d$ and $A$.
\end{theorem}

\pf
It is easy to see that all functionals in $M_0$ are well defined. In particular,
if $\xi$ is the interior vertex $v$ of a pie-shaped triangle $T$ with the curved
edge given by $q(x)=0$ satisfying \eqref{q5}, then the coefficient associated with $v$ is $c_v=s(v)$
regardless of whether we consider $v$ as an element of $D^*_{d-1,T}$ or an
element of $D_{d,\tri_0}$.

For any $\xi\in M_0$ we choose a triangle $T_\xi\in\tri$ such that $\xi\in T_\xi$, 
where we make sure that $T_\xi\in\tri_0$ if $\xi\in D_{d,\tri_0}\setminus\p \Omega$ 
and  $T_\xi\in\tri_P$ if 
$\xi\in \Big(\bigcup _{T\in \tri_{P}}D_{d-1,T}^\ast\Big)\setminus D_{d,\tri_0}$.
Note that for $\xi\in\bigcup _{T\in \tri_{B}} D_{d+1,T}^0$ there is always just one triangle
$T_\xi\in\tri$ that contains $\xi$, and it belongs to $\tri_B$.
If $T_\xi\in\tri_0\cup\tri_B$, then 
$|\lambda_\xi(s)|\le C \|s\|_{L^\infty(T_\xi)}$ for all  $s\in S_{d,0}(\tri)$
by \cite[Theorem~2.6]{LSbook}, where $C$ depends only on $d$. Assume that $T_\xi\in\tri_P$
with the curved side of $T_\xi$ given by $q(x)=0$,
and $s\in S_{d,0}(\tri)$. Then $s|_{T_\xi}=qr$ for some $r\in \mathbb{P}_{d-1}$, and by the
same results of  \cite{LSbook} and Lemma~\ref{stdiv}, 
$|\lambda_\xi(s)|\le C' \|r\|_{L^\infty(T_\xi)}\le C \|s\|_{L^\infty(T_\xi)}$, with a
constant $C$ depending only on $d$ and $B$.
Hence $T_\xi$ is a supporting set of the functional $\lambda_\xi$ satisfying \eqref{K1}
with a constant $K_1$ depending only on $d$ and $B$.

To show that $M_0$ is a minimal determining set and satisfies \eqref{K2} with some 
uniform constant $K_2$, we follow the usual approach (see \cite{LSbook})
of setting the values $c_\xi=\lambda_\xi(s)$, $\xi \in M_0$, 
for any spline $s \in S_{d,0}(\tri)$ and 
showing that the BB-coefficients of $s|_T$ for all $T\in\tri$ can be computed from them consistently and
stably. 

If $T\in\tri_0$, then the BB-coefficients of $s|_T$ are $c_\xi$ for
$\xi\in M_0\cap T=D_{d,T}\setminus\p\Omega$, and zeros for 
$\xi\in D_{d,T}\cap\partial\Omega$, so that no computation is needed. 

If $T\in\tri_P$, then the BB-coefficients of $s|_T$ can be computed by the multiplication
of the BB-form of degree $d-1$ by the BB-form of the quadratic polynomial $q$, see the
explicit formulas \eqref{qmult} below. In view of \eqref{q3} and \eqref{BBstab}, we have
$\|s|_T\|_{L^\infty(T)}\le A \max_{\xi\in D_{d-1,T}^\ast}|\lambda_\xi(s)|$, and 
the BB-coefficients of $s|_T$ are bounded by $\|s|_T\|_{L^\infty(T)}$ times a constant
depending only on $d$. Note that $D_{d-1,T}^\ast\cap(D_{d,\tri_0}\setminus\p \Omega)=\{v\}$,
where $v$ is the interior vertex of $T$. Since $q(v)=1$, we see that the functional 
$\lambda_v(s)=s(v)$ is well-defined regardless of whether we associate it with $v$ as
element of  $D_{d-1,T}^\ast$ or $D_{d,\tri_0}\setminus\p \Omega$.
Hence the BB-coefficient $s(v)$ of $s|_T$ at $v$ is computed consistently.
If $T$ shares an interior edge $e$ with another pie-shaped triangle $T'$, then
Condition (c) of Section~\ref{spacesC0} implies that the curved edges of both $T$ and
$T'$ are given by the same equation $q(x)=0$. Therefor the BB-coefficients of $s$ for
domain points in $D_{d-1,T}^\ast\cap D_{d-1,T'}^\ast$ on the common edge of $T$ and
$T'$ are computed consistently by using the formulas \eqref{qmult} on either $T$ or
$T'$.

For a buffer triangle $T\in\tri_B$ it is easy to see that the coefficients 
$c_\xi=\lambda_\xi(s)$, $\xi \in D_{d+1,T}^0 $ give us the interior part of the 
B\'ezier net of the patch $s\vert_T$, as well as the BB-coefficients on any edges shared
with other buffer triangles. Any BB-coefficients of $s\vert_T$ at the domain points on the boundary of
$\Omega$ are zero, and those on the edges shared with pie-shaped triangles have been computed
already. If the edge $e$ of $T$ is shared with a triangle $T'\in\tri_0$, then $s|_e$ has
already been determined as the restriction of the polynomial $s\vert_{T'}$ of degree $d$.
Therefore we can obtain the BB-coefficients of $s\vert_T$, as a polynomial of degree $d+1$,
for domain points in $e$ by the standard \emph{degree raising} formulas 
\cite[Section 2.15]{LSbook}, which is a stable process. 

Thus, we have shown that  $M_0$ is a minimal determining set. For each $\xi\in M_0$, let $s_\xi$
denote the corresponding dual basis function in $S_{d,0}(\tri)$. By inspecting the above
arguments it is easy to see that 
$$
\supp s_\xi=\cup\{T\in\tri: \xi\in T\}$$
(In particular, \eqref{qmult} shows that
for $\xi$ in the interior of a pie-shaped triangle $T$ the BB-coefficients of $s_\xi|_T$ on 
the edges shared with buffer triangles are zero.) Hence 
$$
\Lambda_T=T\cap M_0\quad\text{for all}\; T\in\tri,$$
which shows that $M_0$ is 1-local according to Definition~\ref{localMDS}, 
with the supporting sets $T_\xi$ described above. It is easy to see that 
$|\Lambda_T|\le{d+2\choose 2}$ if $T\in\tri_0$, $|\Lambda_T|={d+1\choose 2}$ if 
$T\in\tri_P$, and $|\Lambda_T|\le{d+2\choose 2}+d-2$ if $T\in\tri_B$, which implies the
following bound for the covering number:
$$
\kappa_\Lambda\le {d+2\choose 2}+d-2.$$

By inspecting again the above argumentation that $M_0$ is a minimal determining set we
conclude that \eqref{K2} is satisfied with $K_2$ depending only on $d$ and  the constant $A$ 
of \eqref{q3}.
 \eop

The following statement about corresponding basis \eqref{basis} follows immediately.

\begin{corollary}\label{Bbasis}
The dual basis functions $s_\xi:\xi\in M_0$ defined by the condition 
$\lambda_\xi(s_\zeta)=\delta_{\xi\zeta}$ have local support and satisfy the stability estimates of
Proposition~\ref{Lpstab} with constants depending only on $d,A,B$.
\end{corollary}

\section{Implementation of the FEM}\label{implementationC0}

In this section we briefly discuss the implementation aspects of the finite element method
for solving second order elliptic problems on piecewise conic domains using the Bernstein-B\'ezier basis
of the previous section.

\subsection{The conics \label{conics}}

\begin{figure}
\centering
\includegraphics[width=0.5\textwidth]{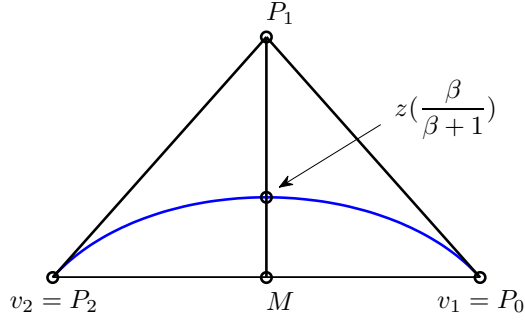} 
\caption{
Conic as a rational B\'ezier curve.
}\label{conicbezier}
\end{figure}

A convenient method to represent the piecewise conic boundary is provided by the rational B\'ezier
curves. 
Given three control points $P_0,{P}_1,{P}_2\in\RR^2$, the 
\emph{quadratic rational B\'ezier curve} can be described by
\begin{equation}\label{Bcurve}
B(t) = \frac{P_0B_0^2 (t) + \beta {P}_1B_1 ^2 (t) + 
{P}_2B_2^2 (t)}{B_0^2 (t) + \beta B_1^2 (t) + B_2^2 (t)},\quad
0 \leq t \leq 1 , 
\end{equation}
where $\beta >0$ is a weight and $B_i ^2 (t) ={{2}\choose{i}} t^i (1-t)^{2-i}$, $i=0,1,2$ 
are quadratic Bernstein polynomials. The curve $B(t)$ goes through $P_0$ and $P_2$ such that 
the tangents at these points are parallel to the segments $\overline{P_0{P}_1}$ and 
$\overline{P_1{P}_2}$, respectively. Moreover, $B(t)$ is contained in the
triangle with vertices $P_0,P_1,P_2$, see Figure~\ref{conicbezier}.
According to \cite[Lemma 4.5]{HoschekBook} 
the curve is a parabolic, elliptic or  hyperbolic arc if 
$\beta =1$, $\beta <1$ or $\beta >1$, respectively.
Let $M$ be the mid-point of the straight line segment $\overline{P_0{P}_2}$. Then the point
$$
(1-s)M+s{P}_1,\quad\text{where }\; s := \tfrac{\beta }{\beta +1},$$
lies on the conic arc \cite{HoschekBook}, %
as shown in Figure~\ref{conicbezier}.

\subsubsection*{Conversion to bivariate Bernstein-B\'ezier form}
\begin{figure}
\centering
\psfrag{c}{$\hspace*{-2pt}\omega$}
\includegraphics[width=0.5\textwidth]{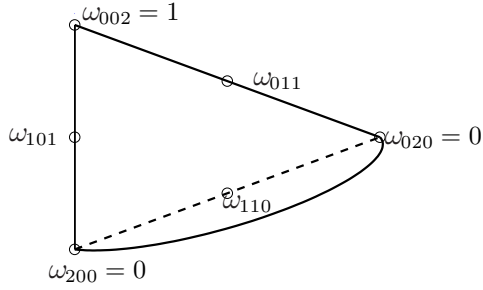} 
\caption{
BB-coefficients of $q$ over a pie-shaped triangle.
}\label{BezierNetq}
\end{figure}

For the sake of computation of the system matrices in the finite element method we need to convert the
conic from the parametric representation \eqref{Bcurve} into implicit form  $q(x)=0$ where $q$ is a 
quadratic polynomial written in BB-form as
$$
q = \sum _{i+j+k=2} \omega_{ijk}B_{ijk} ^2,\quad \omega_{ijk}\in\RR,$$
with respect to the triangle $T^*$ associated with a pie-shaped triangle $T\in\tri_P$, 
see  Figure~\ref{BezierNetq}. Let $v_1=P_0$ and $v_2=P_2$ be the boundary vertices of $T$, and let
$v_3$ be its interior vertex.
Since $q(v_1)=q(v_2)=0$, we have $\omega_{200}=\omega_{020}=0$, and arrive at 
\begin{equation}\label{quadratic1}
q = \omega_{110}B_{110} ^2 + \omega_{101}B_{101} ^2 + \omega_{011}B_{011} ^2 + \omega_{002}B_{002} ^2,
\end{equation}
or, more explicitly,
$$
q = 2(\omega_{110}b_1 b_2  + \omega_{101}b_1 b_3  + \omega_{011}b_2 b_3) + \omega_{002}b_3 ^2.$$

To compute the yet unknown coefficients of $q$, %
we obtain three
linear equations from the facts that the point $z=(1-s)M+s{P}_1$ lies on the conic, and the tangential directional 
derivatives of $q$ at $P_i$, $i=0,1$, vanish. 
Let $\alpha _1 , \alpha _2 , \alpha _3 $ be the 
barycentric coordinates of the control point ${P}_1$ w.r.t.\ $T^*$. Then the barycentric coordinates
of $z$ are
$\tfrac{1-s}{2}+ s\alpha _1 ,\;\tfrac{1-s}{2}+s\alpha _2,\;s\alpha _3$,
so that $q(z)=0$ leads to
the equation
$$
(\tfrac{1-s}{2}+ s\alpha _1)(\tfrac{1-s}{2}+ s\alpha _2)\omega_{110} + 
 s\alpha _3(\tfrac{1-s}{2}+ s\alpha _1)\omega_{101} + s\alpha _3(\tfrac{1-s}{2}+ s\alpha _2)\omega_{011} 
 + \tfrac{1}{2}(s\alpha _3)^2\omega_{002} =0.
$$
It is easy to see that $(\alpha _1 -1 , \alpha _2 , \alpha _3 )$ and $(\alpha _1 , \alpha _2 -1 , \alpha _3 )$ 
are the directional coordinates of the vectors $P_0 {P}_1$ and ${P}_2 {P}_1$ in the terminology of 
\cite[Section~2.6]{LSbook}. Since 
$$
\p_{P_0 {P}_1}q(v_1)=\p_{P_2 {P}_1}q(v_2)=0,$$
we obtain using \cite[Theorem~2.12]{LSbook},
$$
\alpha_2 \omega_{110}+\alpha_3\omega_{101}=0\;\text{ and }\;
\alpha_1 \omega_{110}+\alpha_3\omega_{011}=0,$$
which leads to
$$
\Big(\big(\tfrac{1-s}{2}\big)^2- s^2\alpha _1\alpha _2\Big)\omega_{110} = - \tfrac{1}{2}(s\alpha_3)^2\omega_{002}.
$$
We observe that $\alpha _3 =0$ if and only if ${P}_1$ lies on the straight line $P_0 {P}_2$, which means
that $q$ is zero on this line. This implies that the conditions  $\alpha _3 =0$ and $\omega_{110}=0$
are equivalent to each other and never hold for triangles in $\tri_P$ because their boundary edges are not straight.
Since $0<s<1$, it follows that $\omega_{002}=0$ if and only if 
$\big(\tfrac{1-s}{2}\big)^2- s^2\alpha _1\alpha _2=0$, or using the parameter $\beta$ of \eqref{Bcurve},
$4\beta^2\alpha _1\alpha _2=1$. Since $\omega_{002}=q(v_3)\ne 0$ by our assumptions, we assume without
loss of generality that $q(v_3)=1$, which implies
$$
4\beta^2\alpha _1\alpha _2\ne1$$
and the following formulas for the coefficients of $q$ in \eqref{quadratic1},
\begin{equation}\label{qcurve}
\omega_{110} = -\alpha_3 \mu,\quad \omega_{101}= \alpha_2 \mu, \quad \omega_{011}= \alpha_1 \mu,\quad
\mu:=\frac{2\beta^2\alpha_3}{1-4\beta^2\alpha _1\alpha _2};\quad
\omega_{002}=1.
\end{equation}

Since $T$ is star-shaped w.r.t.~$v_3$, we always have 
$\alpha_1,\alpha_2 > 0$. Moreover, $\alpha _3 > 0\,(<0)$ if ${P}_1$ lies inside $T$ (outside $T$).
We show that the condition that $q(x)$ is positive
for all $x$ in $T$ except of the curved edge is satisfied if and only if $\mu>0$,
that is 
\begin{equation}\label{qpos}
\frac{\alpha_3}{1-4\beta^2\alpha _1\alpha _2}>0.
\end{equation}
Indeed, if $\alpha _3 > 0$, then $T$ is a proper part of $T^*$. Connect $v_3$ by a straight line with 
any point $z$ on the segment $[v_1,v_2]$. Then the univariate quadratic
polynomial obtained by restricting $q$ to this straight line is positive at $v_3$, negative at $z$
and zero at the point of the curved edge it crosses. Hence it cannot be zero at any other point of $T$.
If $\alpha _3 < 0$, then $T$ contains $T^*$ and is convex. By \cite[Theorem 3.6]{LSbook}, $q$ is
positive everywhere in $T^*$ except the vertices $v_1,v_2$ if and only if the coefficients
$\omega_{110},\omega_{101},\omega_{011}$ are all positive, which is equivalent to  $\mu>0$. 
As a conic, the curve $q(x)=0$ cannot have another branch  inside $T\setminus T^*$.

\subsubsection*{Multiplication by $q$ on a pie-shaped triangle}
We now work out formulas for the BB-coefficients $a_{ijk}$ of $s|_T$ in terms of the BB-coefficients of
 $q$ and those of the polynomial $p$ such that $s|_T=qp$. The formulas  \eqref{qmult} have already been used in
the proof of Theorem~\ref{S0dMDS}, and they will again be needed in Section~\ref{assembly}.  

Let
$$
q = \omega_{110}B_{110} ^2 + \omega_{101}B_{101} ^2 + \omega_{011}B_{011} ^2 + B_{002} ^2$$
and 
$$
s|_T=q\sum_{i+j+k=d-1}c_{ijk}B^{d-1}_{ijk}=\sum_{i+j+k=d+1}a_{ijk}B^{d+1}_{ijk}.$$
Since
\begin{align*}
qB^{d-1}_{ijk}
=\;&2\tfrac{(i+1)(j+1)}{d(d+1)}\omega_{110}B^{d+1}_{i+1,j+1,k}
  +2\tfrac{(i+1)(k+1)}{d(d+1)}\omega_{101}B^{d+1}_{i+1,j,k+1}\\
  &+2\tfrac{(j+1)(k+1)}{d(d+1)}\omega_{011}B^{d+1}_{i,j+1,k+1}
	+\tfrac{(k+1)(k+2)}{d(d+1)}B^{d+1}_{i,j,k+2},
\end{align*}
we get the following formula for the BB-coefficients of $s|_T$,
\begin{equation}\label{qmult}
a_{ijk}=\tfrac{2ij\,\omega_{110}}{d(d+1)}c_{i-1,j-1,k}
+\tfrac{2ik\,\omega_{101}}{d(d+1)}c_{i-1,j,k-1}
+\tfrac{2jk\,\omega_{011}}{d(d+1)}c_{i,j-1,k-1}
+\tfrac{(k-1)k}{d(d+1)}c_{i,j,k-2},
\end{equation}
where $c_{r,s,t}:=0$ if at least one of the indices $r,s,t$ is negative.

\subsection{Assembly of the finite element linear system}\label{assembly}
To compute the solution $\tilde u$ of \eqref{discr1} we follow the standard Galerkin scheme and 
expand $\tilde u$ as linear combination of the basis
functions $s_\xi$, $\xi\in M_0$, of Section~\ref{sMDS},
$$
\tilde u=\sum_{\xi\in M_0}c_\xi s_\xi,\quad c_\xi\in\RR.$$
Then  \eqref{discr1} is equivalent to the square linear system
$$
\sum_{\xi\in M_0}c_\xi a(s_\xi,s_\zeta)=(f,s_\zeta)\;\text{ for all }\zeta\in M_0$$
with a positive definite system matrix $[a(s_\xi,s_\zeta)]_{\xi,\zeta\in M_0}$.
For the elliptic equations of second order we have
$$%
a(u,v)=\int_\Omega (\nabla u \cdot A \nabla v +v{\bf b} \cdot \nabla u + cuv)dx,
$$%
where $A$ is a matrix, ${\bf b}$ a vector and $c$ a scalar, each of them in general depending on $x$.
Hence, assembling the matrix requires the computation of the 
following system matrices:  the stiffness matrix
$\mathcal{S}$, the convective matrix $\mathcal{B}$ and the mass matrix $\mathcal{M}$,
whose entries are given by 
\begin{equation*}
\mathcal{S}_{\xi\zeta}=\int _{\Omega} \nabla s_{\xi} \cdot  A\nabla s_{\zeta}dx,\quad 
\mathcal{B}_{\xi\zeta}=\int _{\Omega}  s_{\zeta} {\bf b} \cdot \nabla s_{\xi}dx,\quad
\mathcal{M}_{\xi\zeta}=\int _{\Omega} c s_{\xi} s_{\zeta}dx,\quad \xi,\zeta\in M_0,
\end{equation*}
while the right hand side of the system requires the computation of the load vector $\mathcal{L}$ with entries
$$
\mathcal{L}_{\xi}=\int _{\Omega}fs_{\xi}dx,\quad \xi\in M_0.$$

By integrating over all triangles $T\in\tri$ we can reduce the assembly problem to the computation of 
the \emph{element level system matrices} and \emph{load vector}
\begin{align*}%
&\hat{\mathcal{S}}_T = 
\Big[\int _T \nabla B^{d+i}_{\xi}\cdot  A\nabla B^{d+i}_{\zeta}dx\Big]_{\xi,\zeta\in D_{d+i,T^*}},\;
 \hat{\mathcal{B}}_T = 
 \Big[\int _T  B^{d+i}_{\xi}{\bf b}\cdot \nabla B^{d+i}_{\zeta}dx\Big]_{\xi,\zeta\in D_{d+i,T^*}},\\
&\hat{\mathcal{M}}_T = \Big[\int _T c B^{d+i}_{\xi} B^{d+i}_{\zeta}dx\Big]_{\xi,\zeta\in D_{d+i,T^*}},
 \;\hat{\mathcal{L}}_{T}=\Big[\int _{T}fB^{d+i}_{\xi}dx\Big]_{\xi\in D_{d+i,T^*}},
\end{align*}
where $i=0$ if $T\in\tri_0$, $i=1$ if
$T\in\tri_P\cup\tri_B$, and $T^*$ is defined as above for $T\in\tri_P$ and simply as
$T^*=T$ otherwise. Indeed, it is easy to see that
\begin{equation}\label{transf}
\mathcal{S} = \mathcal{T}^t\hat{\mathcal{S}}\mathcal{T},\; 
 \mathcal{B} = \mathcal{T}^t\hat{\mathcal{B}}\mathcal{T},\;
 \mathcal{M} = \mathcal{T}^t\hat{\mathcal{M}}\mathcal{T},\; 
 \mathcal{L} = \mathcal{T}^t\hat{\mathcal{L}},
\end{equation}
where
\begin{align*}
\hat{\mathcal{S}}&:=\diag( \hat{\mathcal{S}}_{T},T\in\tri),\quad
\hat{\mathcal{B}}:=\diag( \hat{\mathcal{B}}_{T},T\in\tri),\\
\hat{\mathcal{M}}&:=\diag( \hat{\mathcal{M}}_{T},T\in\tri),\quad
\hat{\mathcal{L}}:=[ \hat{\mathcal{L}}_{T}^t,T\in\tri],
\end{align*}
and $\mathcal{T}$ is the \emph{transformation matrix} whose columns correspond to the basis functions
$s_\xi$, $\xi\in M_0$, and consist of the BB-coefficients of $s_\xi|_T$ on all triangles $T\in\tri$, and 
$\mathcal{T}^t$ is the transpose of $\mathcal{T}$. 

The structure of the transformation matrix is rather simple. Let $\tau_\xi$ denote the column of
$\mathcal{T}$ corresponding to $\xi\in M_0$, such that $\mathcal{T}=[\tau_\xi]_{\xi\in M_0}$.
If $\xi\in D_{d,\tri_0}$ and $\xi$ is not contained in any pie-shaped or buffer triangle, then
the only non-zero entries of $\tau_\xi$ are the ones at every row corresponding to $\xi\in D_{d,T}$ 
for some $T\in\tri_0$. Similarly, if $\xi\in D^0_{d+1,T}$ for some $T\in \tri_B$, then 
$\tau_\xi$  also consists of ones and zeros, with ones at all rows corresponding to the same 
$\xi$ in one or more triangles of $\tri_B$. 
If $\xi=\xi_{ijk}\in D^*_{d-1,T}$ for some $T\in \tri_P$, then  four entries of $\tau_\xi$ 
will be non-zero in rows corresponding to $\xi_{i+1,j+1,k}$, $\xi_{i+1,j,k+1}$, $\xi_{i,j+1,k+1}$, 
 $\xi_{i,j,k+2}$ in $D_{d+1,T^*}$. These entries can be determined from the
equation \eqref{qmult} as 
$$
\tfrac{2(i+1)(j+1)\,\omega_{110}}{d(d+1)},\quad\tfrac{2(i+1)(k+1)\,\omega_{101}}{d(d+1)},\quad
\tfrac{2(j+1)(k+1)\,\omega_{011}}{d(d+1)},\quad\tfrac{(k+1)(k+2)}{d(d+1)},$$ 
respectively. If $i=0$ and $j\ne0$, then the domain points 
$\xi_{i,j+1,k+1}=\xi_{0,j+1,k+1}$ and  $\xi_{i,j,k+2}=\xi_{0,j,k+2}$ lie on 
an edge of $T$ shared with a triangle $T'$ in $\tri_P\cup\tri_B$, and hence the two
entries of $\tau_\xi$ corresponding to these points in $T'$ will be filled with the
numbers 
$\tfrac{2(j+1)(k+1)\,\omega_{011}}{d(d+1)}$, $\tfrac{(k+1)(k+2)}{d(d+1)}$,
respectively. Similarly, there are two additional nonzero entries in the case 
$i\ne0$, $j=0$.
If $i=j=0$, then $k=d-1$ and  $\xi_{i,j,k+2}=\xi_{0,0,d+1}$ is the interior vertex $v$
of $T$. It is clear from the above that the nonzero entries of $\tau_v$ are:
(a) the ones at all rows
corresponding to $v$ in any types of triangles containing this vertex, 
(b) the numbers 
$\tfrac{2\,\omega_{110}}{d(d+1)}$, $\tfrac{2\,\omega_{101}}{d+1}$, 
$\tfrac{2\,\omega_{011}}{d+1}$  
placed in the rows corresponding to $\xi_{1,1,d-1}$, $\xi_{1,0,d}$, $\xi_{0,1,d}$, 
respectively, in all pie-shaped triangles attached to $v$ 
(where $\omega_{110},\omega_{101},\omega_{011}$ depend on the pie-shaped triangle), 
and (c) the number
$\tfrac{2\,\omega_{101}}{d+1}$ or $\tfrac{2\,\omega_{011}}{d+1}$
in the row corresponding to $\xi_{1,0,d}$ or $\xi_{0,1,d}$ if it belong to a buffer
triangle attached to the respective pie-shaped triangle.
Finally, assume that  $\xi\in D_{d,\tri_0}\cap T$ for a buffer triangle $T$, but 
$\xi$ is not a vertex of any pie-shaped triangle. Then $\xi=\xi_{ijk}$ for some 
$i,j,k$ with $i+j+k=d$, as element of $D_{d,T}$. (Note that at least one of the indices 
$i,j,k$ is zero because $\xi$ lies on an edge of $T$.) Then by using the well known
\emph{degree raising} formulas \cite[Theorem 2.39]{LSbook}, we see that the non-zero
entries of $\tau_\xi$ are, in addition to ones in the rows corresponding to $\xi$ in
all triangles of $\tri_0$ containing $\xi$, the numbers 
$\tfrac{i+1}{d+1}$, $\tfrac{j+1}{d+1}$, $\tfrac{k+1}{d+1}$
in the rows corresponding to $\xi_{i+1,j,k}$, $\xi_{i,j+1,k}$ and 
$\xi_{i,j,k+1}$  as elements of $D_{d+1,T}$, for each buffer triangle
$T$ containing $\xi$.

Note that the global-local transformation \eqref{transf} can be performed without explicit evaluation
and storage of the transformation matrix $\mathcal{T}$. Instead, suitable routines for the matrix-vector products 
$\mathcal{T}x$ and  $\mathcal{T}^tx$ need to be implemented, and this can be done efficiently with
$\mathcal{O}(d^2)$ computation cost per triangle $T\in\tri$, similar to the 
algorithms in \cite[Section 8.1]{ADS}. 

As shown in \cite{AAD11}, the element level system matrices can be computed with optimal 
$\mathcal{O}(1)$ cost per entry on the non-curved triangles. Indeed, 
the product formula for Bernstein polynomials
$$
B_{ijk}^{d}B_{rst}^{q}=
\frac{{{i+r}\choose{i}}{{j+s}\choose{j}}{{k+t}\choose{k}}}{{{d+q}\choose{d}}}B_{i+r,j+s,k+t}^{d+q}$$
helps to reduce this task to the computation of the \emph{Bernstein-B\'ezier (BB-) moments}
$$
\mu^{T,m}_\xi(g):=\int_{T}g(x)B^{m}_{\xi}(x)\,dx,\qquad \xi\in D_{m,T},$$
where $g$ is one of the functions $A,{\bf b},c$ and $m$ is a number between $2d-2$ and $2d+2$. The
components of the load vector $\hat{\mathcal{L}}_{T}$ are just the moments of $f$ of 
degree $m=d$ or $d+1$. Using \cite[Algorithm 3]{AAD11}, the moment vector 
$\mu^{T,m}(g):=[\mu^{T,m}_\xi(g)]_{\xi\in D_{m,T}}$ can be evaluated with $\mathcal{O}(m^3)$ floating point operations
with the help of an $\mathcal{O}(m^2)$-point Stroud quadrature, which delivers sufficient accuracy
for the finite element approximation. 

If $T\in\tri_P$, then the same approach can be used, where we only need to figure out how to compute
the BB-moments.

\subsection{BB-moments on curved triangles}\label{QuadraturePieShaped}
We show that for any pie-shaped triangle $T\in\tri_P$ the moment vector  $\mu^{T,m}(f)$ can be computed
with $\mathcal{O}(q^2m)$ %
operations using a quadrature rule with $q^2$ centers, which is
exact for all $f\in \mathbb{P}_{2q-m-1}$, $q> m/2$. This means that $q=\mathcal{O}(d)$ points 
guarantee sufficient accuracy, and hence the element level system matrices are computed on pie-shaped
triangles also with optimal computation cost $\mathcal{O}(1)$ per entry  \cite{AAD11}.
 
As before let $v_1,v_2$ be the boundary vertices of $T$, and $v_3$ its interior vertex. We consider
a version of Duffy transformation $\Phi$ \cite[Section 3.1]{AAD11} that maps any point $(t_1,t_2)\in\RR^2$ into 
$\Phi(t_1,t_2)=b_1v_1+b_2v_2+b_3v_3$, where
\begin{equation}\label{Duffy}
b_1=(1-t_1)t_2,\quad b_2=t_1t_2,\quad b_3=1-t_2.
\end{equation}
It is easily seen that the straight triangle $T^*$ with vertices $v_1,v_2,v_3$ is the image of the unit
square, that is $T^*=\Phi([0,1]^2)$, with $v_1=\Phi(0,1)$, $v_2=\Phi(1,1)$, $v_3=\Phi(1,0)=\Phi(0,0)$. 
For any constant  $\tau$ the
image of the straight line $t_1=\tau$ is the straight line going 
through $v_3$ and the point $(1-\tau)v_1+\tau v_2$ lying on the line through $v_1,v_2$. The pie-shaped
triangle $T$ is the image of 
$$
\hat T:=\{(t_1,t_2)\in \RR^2:0\le t_1\le 1,\;0\le t_2\le \phi(t_1)\},$$
where $\phi$ is a continuous function which is well defined because $T$ is star-shaped with respect to $v_3$.
Assuming that the curved edge of $T$ is given by the equation $q(x)=0$, with $q$ defined by
\eqref{quadratic1}, $\phi(\tau)$ is easy to compute for any $0\le\tau\le1$ as the first positive root of the univariate quadratic 
polynomial
\begin{align*}
\tilde q(t)&=q(v_3+t((1-\tau)v_1+\tau v_2 -v_3))\\
&=\omega_{002}B^2_0(t)+\big(\omega_{101}(1-\tau)+ \omega_{011}\tau\big)B^2_1(t)
+2\omega_{110}(1-\tau)\tau B^2_2(t),
\end{align*}
where we use quadratic Bernstein polynomials $B^2_i$ as in \eqref{Bcurve}.

Since the Jacobian of $\Phi$ is $2|T^*|t_2$, where $|M|$ denotes the area of a set $M\subset \RR^2$,
we can compute any integral over $T$ as
$$
\int_Tf(x)\,dx=2|T^*|\int_{0}^1dt_1\int_0^{\phi(t_1)}f(\Phi(t_1,t_2))t_2\,dt_2.$$
By applying Gauss-Legendre, respectively Gauss-Jacobi quadrature with $q$ points,
$$
\int_0^1g(s)\,ds\approx\sum_{j=1}^q w_jg(\xi_j),\qquad 
\int_0^1sg(s)\,ds\approx\sum_{j=1}^q \tilde w_jg(\tilde \xi_j),$$
to the  integrals in $t_1$, respectively $t_2$, we obtain a positive quadrature rule with $q^2$ centers
\begin{equation}\label{gStroud}
\int_Tf(x)\,dx\approx 2|T^*|\sum_{\mu=1}^q w_\mu\phi^2(\xi_\mu)\sum_{\nu=1}^q \tilde w_\nu
f(\xi_{\nu\mu}),\quad\text{with }\; \xi_{\nu\mu}:=\Phi(\xi_\mu,\phi(\xi_\mu)\tilde \xi_\nu),
\end{equation}
which is exact whenever $f\in\mathbb{P}_{2q-1}$ since $f(\Phi(t_1,t_2))$ is a univariate polynomial of
degree $2q-1$ in each of $t_1,t_2$ in that case.

Similar to \cite[Lemma 1]{AAD11} the bivariate Bernstein polynomials with respect to $T^*$ 
can be factorized as
$$
B^m_{ijk}(x)=B^{m-k}_j(t_1)B^m_k(t_2),\quad x=\Phi(t_1,t_2),$$   
where $B^n_\nu(t)={n\choose \nu}t^\nu(1-t)^{n-\nu}$, $0\le\nu\le n$, are the univariate Bernstein polynomials.
Hence, by using \eqref{gStroud} we obtain the following approximation of the moments for all
$i+j+k=m$,
\begin{equation}\label{momquad}
\mu^{T,m}_{ijk}(f)\approx 2|T^*|\sum_{\mu=1}^q w_\mu\phi^2(\xi_\mu)B^{m-k}_j(\xi_\mu)
\sum_{\nu=1}^q \tilde w_\nu B^m_k(\phi(\xi_\mu)\tilde \xi_\nu)
f(\xi_{\nu\mu}).
\end{equation}
Since $B^m_{ijk}$ is a polynomial of degree $m$ in $t_2$ and at most $m$ in $t_1$, the formula 
\eqref{momquad} is exact for all $f\in\mathbb{P}_{2q-m-1}$ as long as $q>m/2$, and its structure allows evaluation by
sum factorization as in \cite[Section 3.3]{AAD11}. That is, first the sums
$$
\sigma_{\mu,k}=\sum_{\nu=1}^q \tilde w_\nu B^m_k(\phi(\xi_\mu)\tilde \xi_\nu)f(\xi_{\nu\mu})$$ 
are computed for all $\mu=1,\ldots,q$ and $k=0,\ldots,m$ with $\mathcal{O}(q^2m)$ cost, and then the
numbers
$$
2|T^*|\sum_{\mu=1}^q w_\mu\phi^2(\xi_\mu)B^{m-k}_j(\xi_\mu)\sigma_{\mu,k},
\quad k=0,\ldots,m,\; j=0,\ldots,m-k,$$
give the moments $\mu^{T,m}_{ijk}(f)$, $i=m-j-k$, with the cost of $\mathcal{O}(qm^2)$.

\section{Numerical experiments}\label{examplesC0}
To test the numerical performance of our method we implemented it in  MATLAB and
present in this section three examples 
including the membrane eigenvalue problem and Poisson problem  on different domains with curved boundaries. 
We consider both $h$- and $p$-refinements and compare our results to the 
state-of-the-art software COMSOL Multiphysics. We use version 4.2a of
COMSOL which provides the options to employ the standard isoparametric  elements  of degree up 
to 5 on curved domains in 2D. While the domains in Examples 1 and 2 are smooth (an ellipse and a circle),
we consider a $C^0$ domain bounded by linear and conic pieces in Example 3. 
The numerics confirm in particular the theoretical rate of convergence given 
in Theorem~\ref{bound1}. %
Note that for simplicity we used an easier implementation which does not achieve optimal cost assembly described in
Section~\ref{implementationC0}. (Further details can be found in \cite{Abid_thesis}.)

\subsection*{Example 1: Poisson problem on an ellipse}
Let $\Omega$ be the domain in  $\RR^2$ bounded by the ellipse $x_1^2+6.25x_2^2=1$. Consider the
boundary value problem
\begin{equation}\label{poisson}
\Delta u =f \mbox{ in } \Omega, \quad
u = 0 \mbox{ on } \partial \Omega,
\end{equation}
where $f$ is chosen such that the exact solution of the problem is
$u=e^{0.5(x_1^2+6.25x_2^2)}-e^{0.5}$.

\begin{figure}
\centering
\begin{tabular}{c}
\includegraphics[width=.6\textwidth]{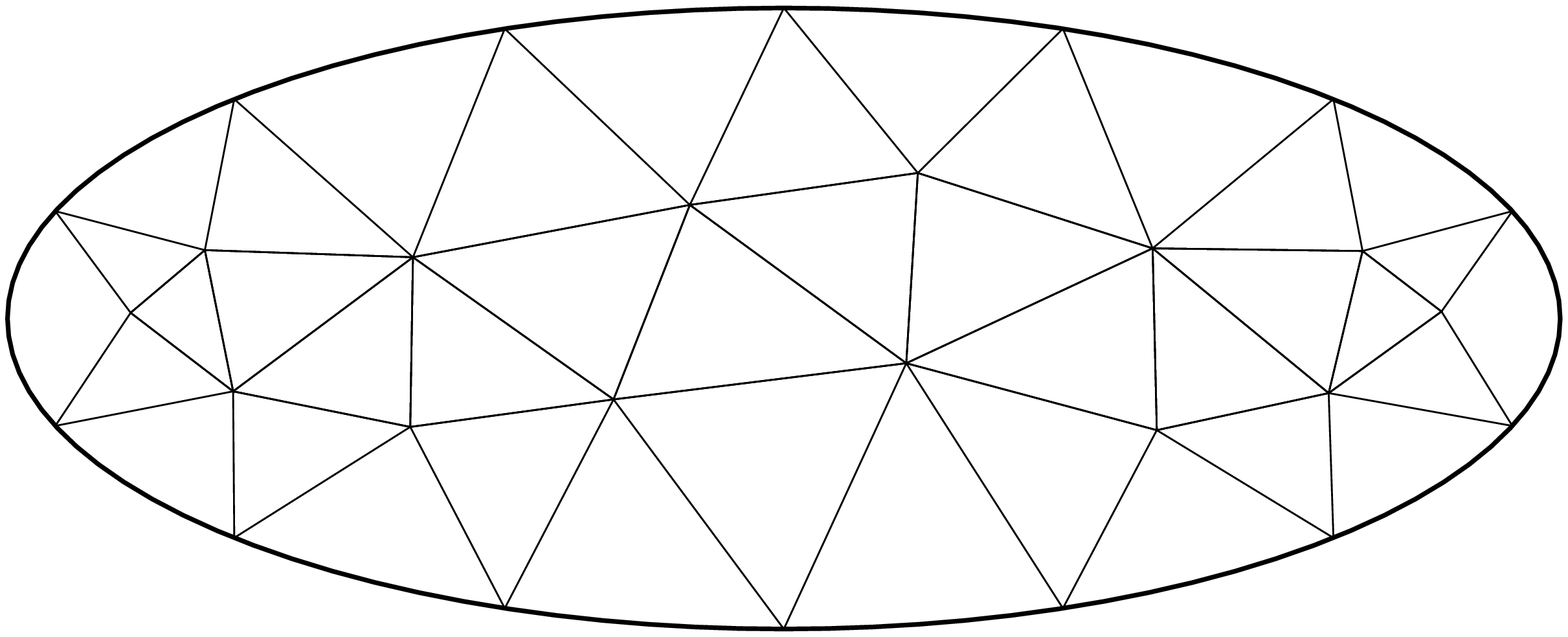} 
\end{tabular}
\caption{
Example 1: Initial triangulation.
}\label{Ex1init}
\end{figure}

\begin{figure}
\centering
\begin{tabular}{c}
\includegraphics[width=.5\textwidth]{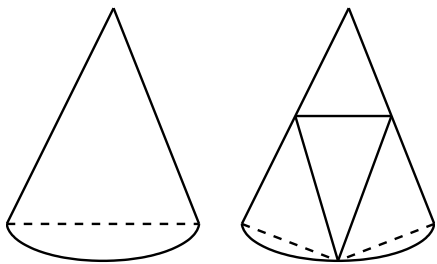} 
\end{tabular}
\caption{
Refinement of a pie-shaped triangle.
}\label{RefinePieShape}
\end{figure}

To ensure a fair comparison of the numerical results 
with COMSOL, we use the initial triangulation shown in Figure~\ref{Ex1init} generated by 
COMSOL. We obtain a sequence of triangulations  by uniform refinement whereby each triangle is subdivided into four
subtriangles by joining the midpoints of all edges. For a pie-shaped triangle 
 we take the midpoint of the curved boundary edge, see Figure~\ref{RefinePieShape}.

\begin{figure}
\centering
\begin{tabular}{c}
\includegraphics[height=0.6\textwidth]{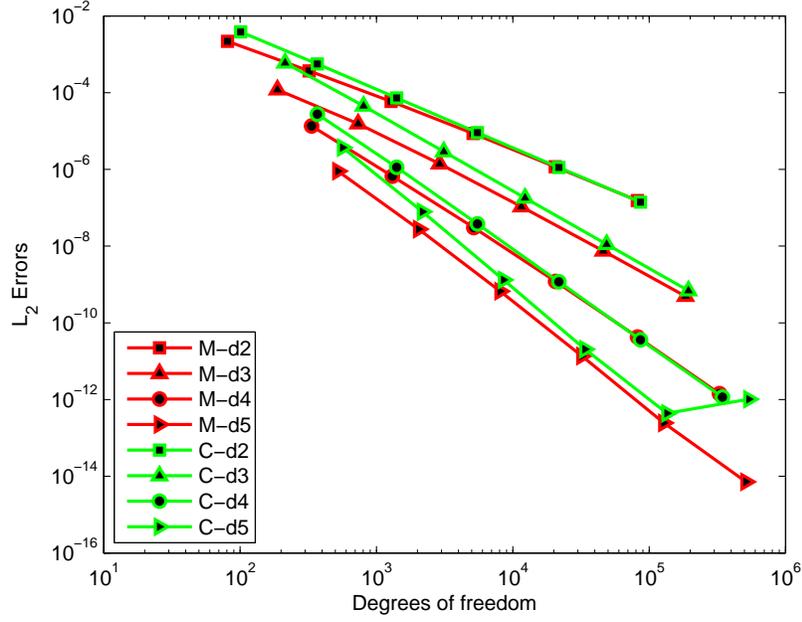} 
\end{tabular}
\caption{
$L_2$ errors in Example 1 using our method (red) and COMSOL (green). The curves indicated by {\tt
d2, d3, d4, d5} correspond to the degrees $d=2,\ldots,5$.
}\label{L2Test2}
\end{figure}

\begin{figure}
\centering
\begin{tabular}{c}
\includegraphics[height=0.6\textwidth]{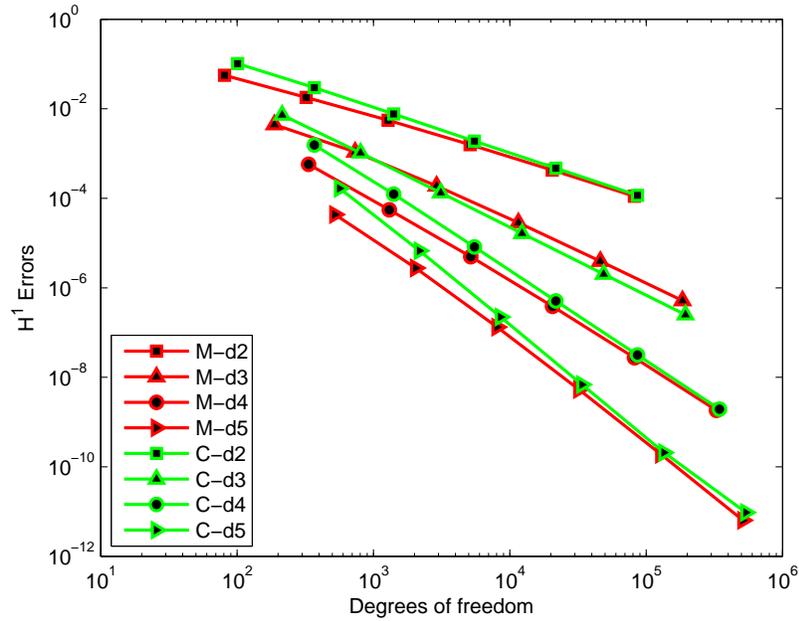} 
\end{tabular}
\caption{
$H^1$ errors in Example 1 using our method (red) and COMSOL (green).
}\label{H1Test2}
\end{figure}

\begin{figure}
\centering
\begin{tabular}{c}
\includegraphics[height=0.6\textwidth]{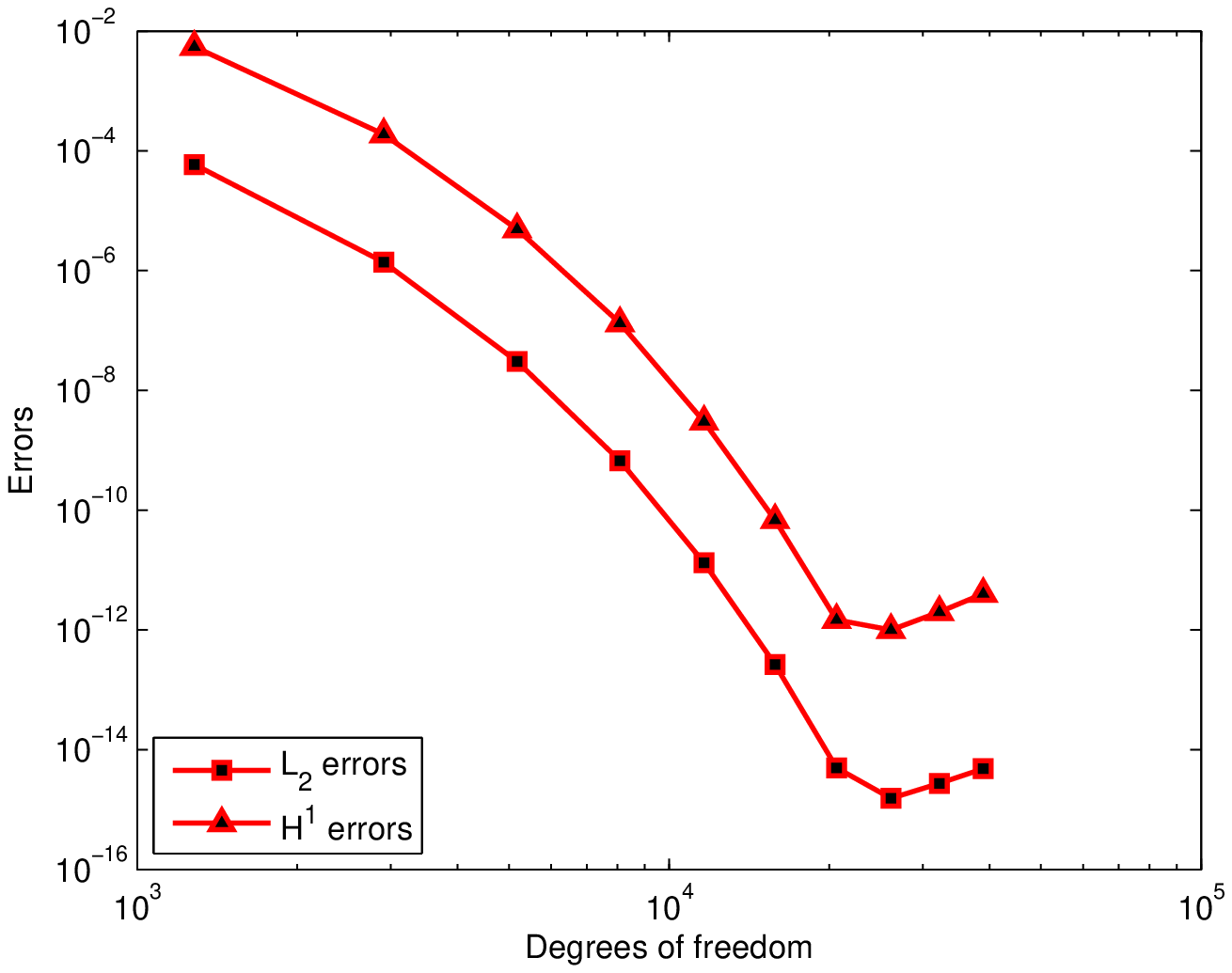} 
\end{tabular}
\caption{
$L_2$ and $H^1$ errors in Example 1 using $p$-refinement.
}\label{pMethodTest2}
\end{figure}

Error plots for
 $\| u-u_N \| _{L_2(\Omega)}$
 and 
 $\| u-u_N \| _{H^1(\Omega)}$ as functions of the number of degrees of freedom $N=\dim S_{d,0}(\tri)$
	are depicted in Figures~\ref{L2Test2} and  
 \ref{H1Test2}, respectively, for $d=2,3,4,5$, 
 where $u_N\in S_{d,0}(\tri)$ is the approximate solution to the problem, and different points on the
 curves correspond to subsequent refinements of the triangulation. 
The results show that both methods behave similarly for all different orders and 
confirm the theoretical estimates of Theorem~\ref{bound1}.
Furthermore, Figure~\ref{pMethodTest2} demonstrates the expected exponential decay of errors of our
method in an
experiment of a ``$p$-refinement'' type \cite{BnS2}, where the triangulation is fixed (as obtained by two subsequent
uniform refinements of the initial triangulation), and the degree $d$ is increased instead, 
increasing this way the number of degrees of freedom.

\subsection*{Example 2: Circular membrane eigenvalue problem}

The free vibrations of a homogeneous membrane are governed by the equation 
\begin{equation}\label{membrane3}
\Delta u + \lambda u =0, \quad x \in \Omega.
\end{equation}
If the membrane is fixed along its boundary then the boundary condition is 
\begin{equation}\label{membr3}
 u =0, \quad x \in \p \Omega,
\end{equation}
which comprises a problem of finding eigenvalues and eigenfunctions of the Laplacian 
under homogeneous  Dirichlet boundary conditions.

We choose $\Omega \subset \RR^2$ to be a unit disk and approximate 15
smallest  eigenvalues for the circular membrane. The exact solution 
to this problem is known \cite{KS84}, with the eigenvalues  given by 
$$\lambda _{m,n}=(j_{m,n})^2, \quad m=0,1,\hdots, \quad n=1,2,\hdots,$$
where $j_{m,n}$ is the $n$-th root of the $m$-th Bessel function $J_m$ of the first kind.

The weak variational formulation corresponding to \eqref{membrane3} and \eqref{membr3}, for $\lambda\in\RR$ 
and $u\ne0$, is
\begin{equation}
\label{membr4}
\mbox{ Find }\,u\in H^1_0(\Omega),\;
\int_\Omega\nabla u\cdot\nabla v=\lambda\int_\Omega uv,\; 
\hbox{$\forall $ $v\in H^1_0(\Omega)$.}
\end{equation}
We discretize this problem in the space $S_{d,0}(\tri)$ for a suitable triangulation of $\Omega$: 
Find $\tilde u\in S_{d,0}(\tri)$,  $\tilde u\ne0$, such that
\begin{equation}
\label{membr5}
\int_\Omega\nabla \tilde u\cdot\nabla \tilde v=\lambda\int_\Omega \tilde u\tilde v,\quad 
\hbox{$\forall \tilde v\in S_{d,0}(\tri)$.}
\end{equation}
Hence if $\left\lbrace s_1,\hdots,s_N\right\rbrace $ is a  for 
$S_{d,0}(\tri)$ according to Theorem~\ref{S0dMDS}, then  
\eqref{membr5} boils down to the matrix equation of the form
 $$\mathcal{S}=\lambda \mathcal{M},$$
 where $\mathcal{S}$ and $\mathcal{M}$ are the stiffness and mass matrices. We solve this generalized
 matrix  eigenvalue problem using MATLAB's built-in command {\tt eig}.

\begin{figure}
\centering
\begin{tabular}{c}
\includegraphics[width=0.4\textwidth]{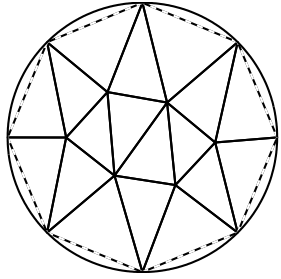} 
\end{tabular}
\caption{
Initial triangulation  of the unit circle for Example 2.
}\label{initialmeshcircle}
\end{figure}

 We follow the same procedure to get a sequence of 
meshes as in Example 1, starting with the initial mesh shown in Figure~\ref{initialmeshcircle}
imported from COMSOL. 
Note that this triangulation does not satisfy \eqref{q1}.

\begin{figure}
\centering
\begin{tabular}{c}
\includegraphics[height=0.6\textwidth]{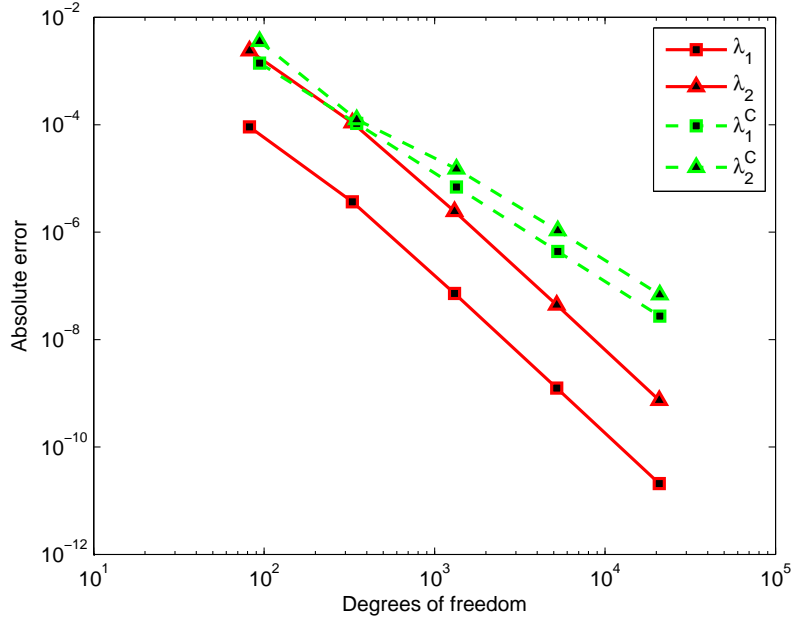} 
\end{tabular}
\caption{
Example 2: Absolute errors for 1st and 2nd eigenvalues by our method $(\lambda _1,\lambda _2)$ and by 
COMSOL ($\lambda ^C_1$,$\lambda ^C_2$) for $d=3$.}\label{eigen1-2-d3}
\end{figure}

\begin{figure}
\centering
\begin{tabular}{c}
\includegraphics[height=0.6\textwidth]{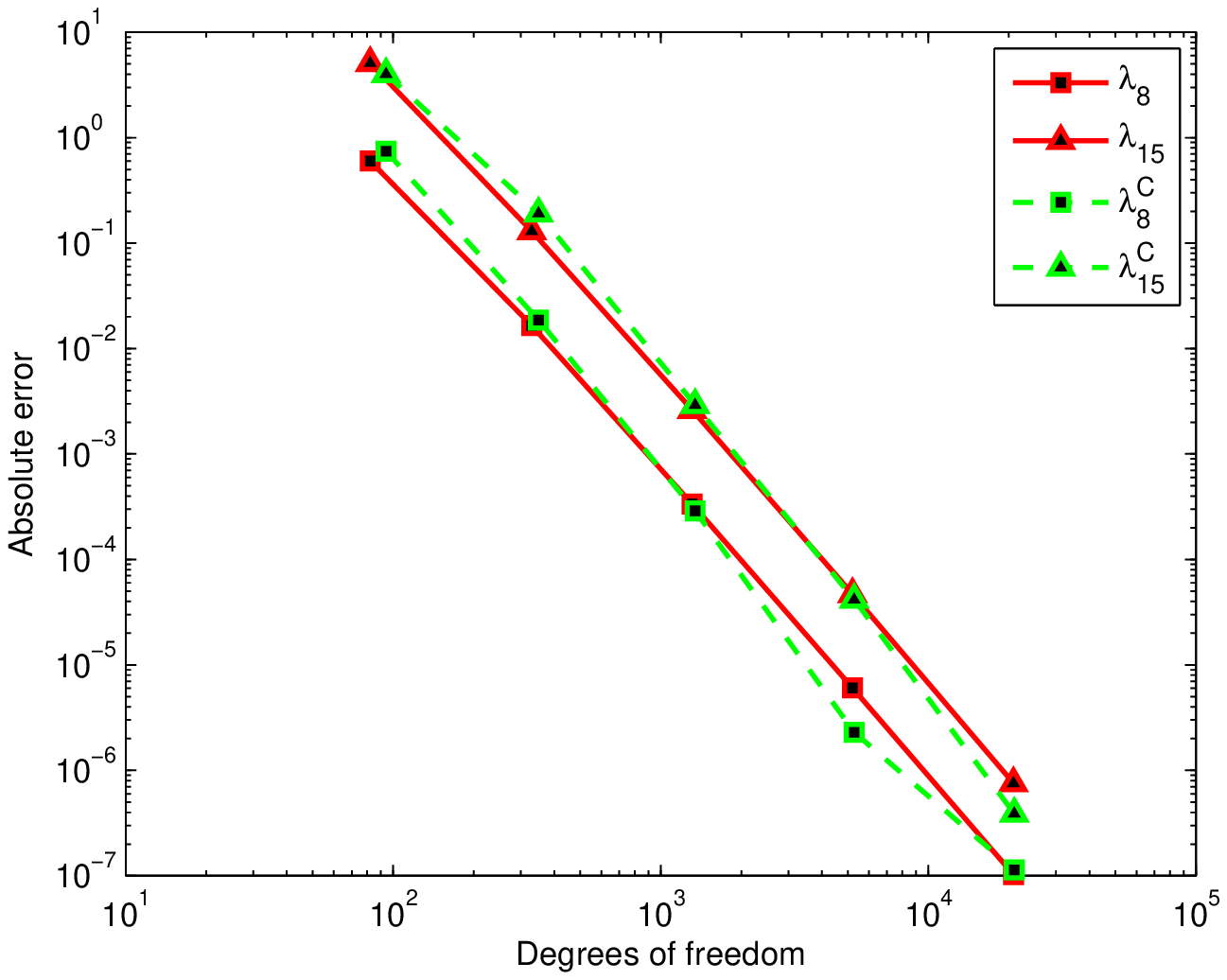} 
\end{tabular}
\caption{Example 2: Absolute errors 
for 8th and 15th eigenvalues for $d=3$ ($\lambda _8,\lambda _{15}$: our method,
$\lambda^C _8,\lambda^C _{15}$: COMSOL).
}\label{eigen8-15-d3}
\end{figure}

\begin{figure}
\centering
\begin{tabular}{c}
\includegraphics[height=0.6\textwidth]{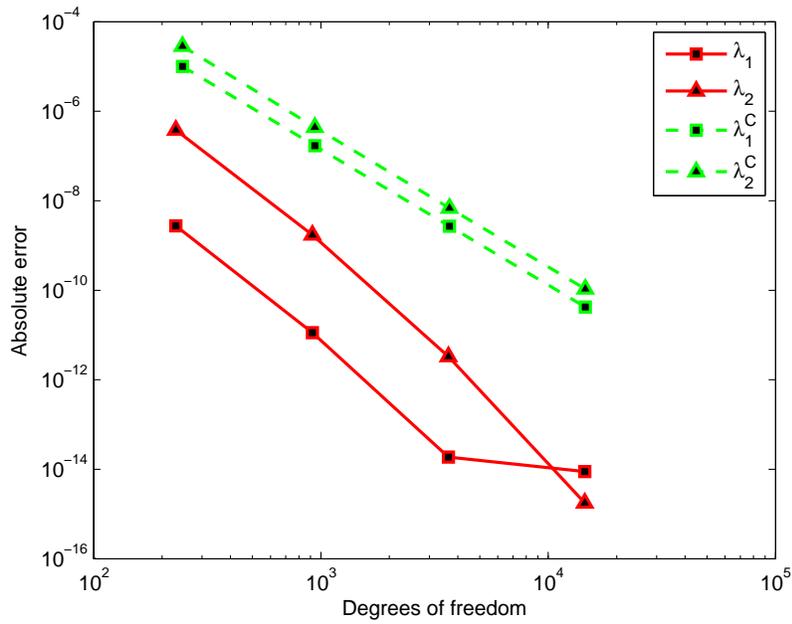} 
\end{tabular}
\caption{Example 2: Absolute errors 
for 1st and 2nd eigenvalues for $d=5$ ($\lambda _1,\lambda _{2}$: our method,
$\lambda^C _1,\lambda^C _{2}$: COMSOL).
}\label{eigen1-2-d5}
\end{figure}

\begin{figure}
\centering
\begin{tabular}{c}
\includegraphics[height=0.6\textwidth]{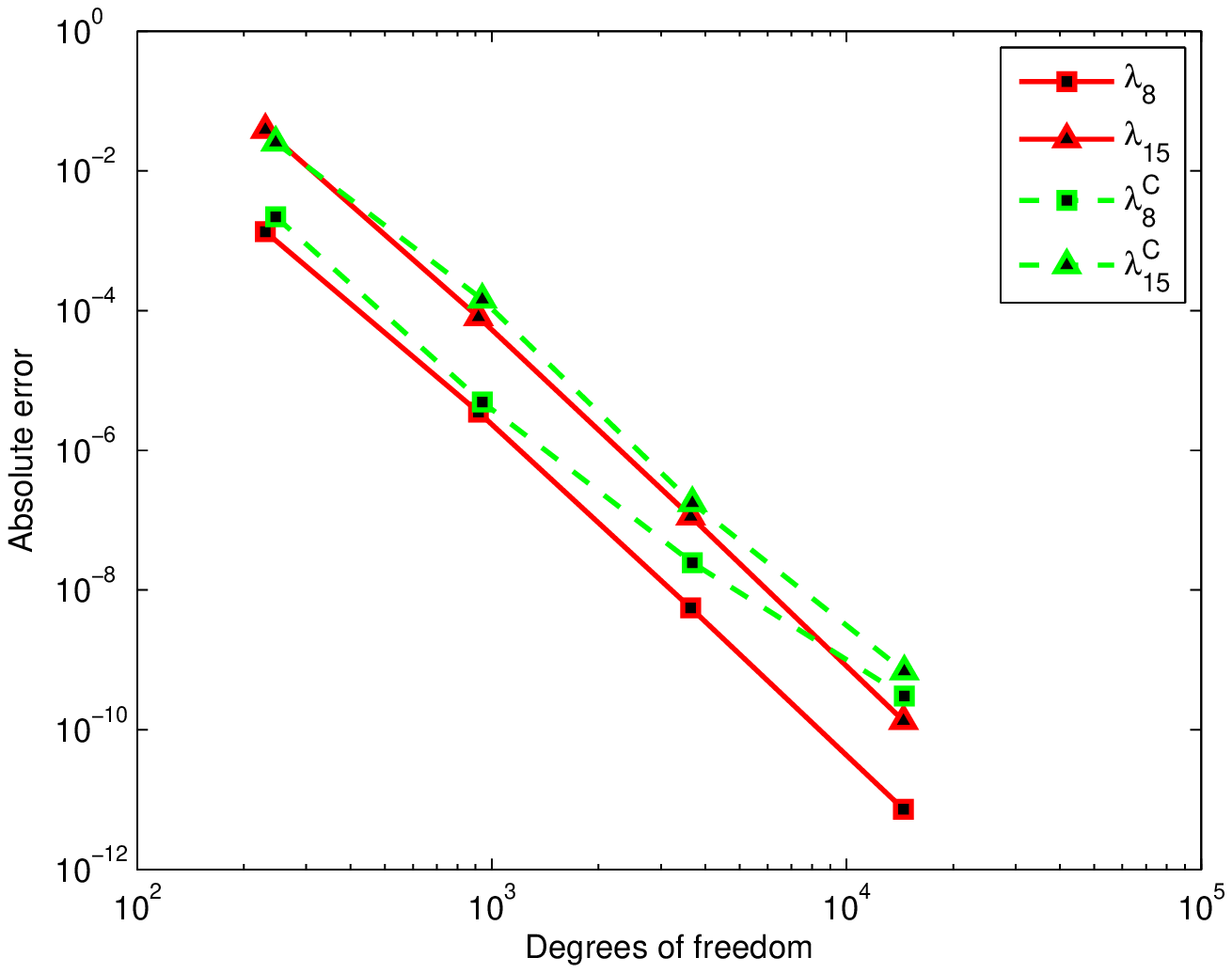} 
\end{tabular}
\caption{Example 2: Absolute errors 
for 8th and 15th eigenvalues for $d=5$ ($\lambda _8,\lambda _{15}$: our method,
$\lambda^C _8,\lambda^C _{15}$: COMSOL).
}\label{eigen8-15-d5}
\end{figure}

\begin{figure}
\centering
\begin{tabular}{c}
\includegraphics[height=0.637\textwidth]{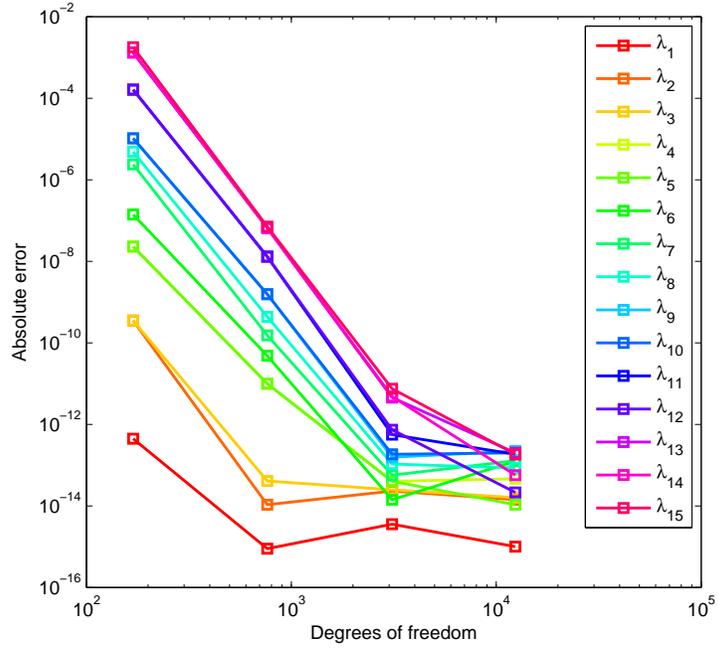} 
\end{tabular}
\caption{
Example 2: Absolute errors for the first 15 eigenvalues using our method  for $d=9$.
}\label{eigenAll-d9}
\end{figure}

In Figures~\ref{eigen1-2-d3}--\ref{eigen8-15-d5} we plot absolute errors for approximating 
various eigenvalues using our implementation and using COMSOL for degree $d=3$ and $d=5$. 
It can be seen that, comparing to COMSOL, our method approximates the first two eigenvalues 
significantly better. 
For the 8th and 15th eigenvalues the results are comparable for $d=3$ 
 but for the higher degree $d=5$ the accuracy achieved using our method is again better. 
Figure~\ref{eigenAll-d9} depicts the errors of our method for the first 
15 eigenvalues for degree $d=9$. 
It is easy to check that the slopes of the curves in Figures~\ref{eigen1-2-d3}--\ref{eigenAll-d9} are
consistent with the expected convergence rate $h^{2d}$ which follows from Theorem~\ref{approx} in view
of \cite[Theorem~3.1]{BO87}.

\begin{figure}
\centering
\begin{tabular}{c}
\includegraphics[height=0.6\textwidth]{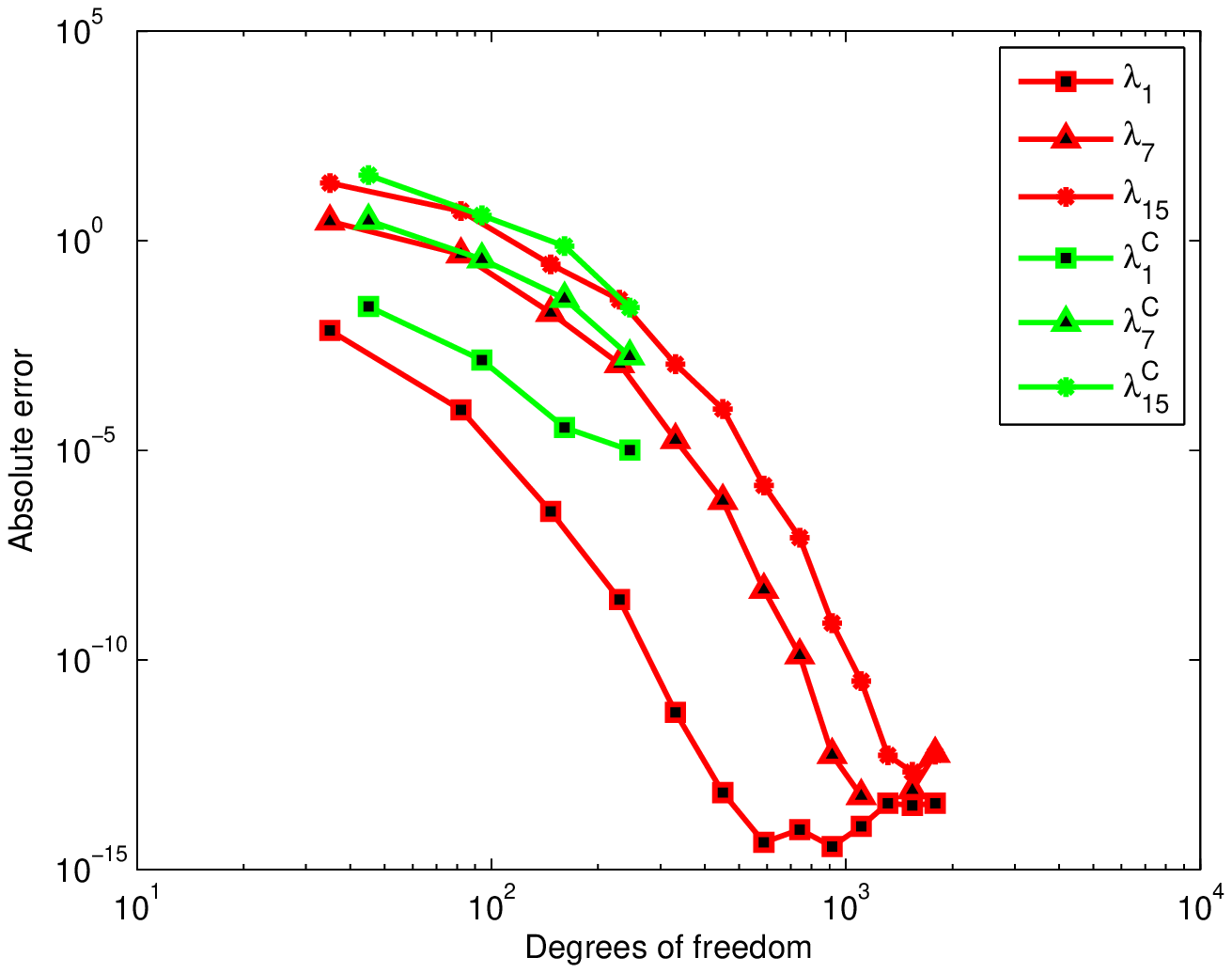} 
\end{tabular}
\caption{
Example 2: Absolute errors for the 1st, 7th and 15th eigenvalues using $p$-refinement
 over initial triangulation shown in Figure~\ref{initialmeshcircle} 
 ($\lambda _1,\lambda _{7},\lambda _{15}$: our method,
$\lambda^C _1,\lambda^C _7,\lambda^C _{15}$: COMSOL).
}\label{threeEig-pref}
\end{figure}

We have also tested the $p$-refinement for this problem on the initial triangulation shown in 
Figure~\ref{initialmeshcircle}. We plot the absolute errors for the 1st, 7th and 15th eigenvalues in 
Figure~\ref{threeEig-pref}. For comparison, the figure also includes the errors obtained with COMSOL up the highest degree 5 it
allows. The rate of convergence is exponential in this case as expected. 

\subsection*{Example 3: Poisson problem on a $C^0$ piecewise conic domain}

Here we consider a domain $\Omega$ with $C^0$ boundary bounded by linear and quadratic boundary segments. 
The boundary is defined by two straight lines 
$x_2 = \pm 2$ and two parabolas $x_1 = \pm (x_2^2 -6)$. 
We consider the Poisson problem \eqref{poisson}, where $f$ is chosen such that
the exact solution is $u=(x_2^2-4)(x_1^2-(x_2^2-6)^2)/100$.

\begin{figure}
\centering
\begin{tabular}{c}
\includegraphics[height=0.6\textwidth]{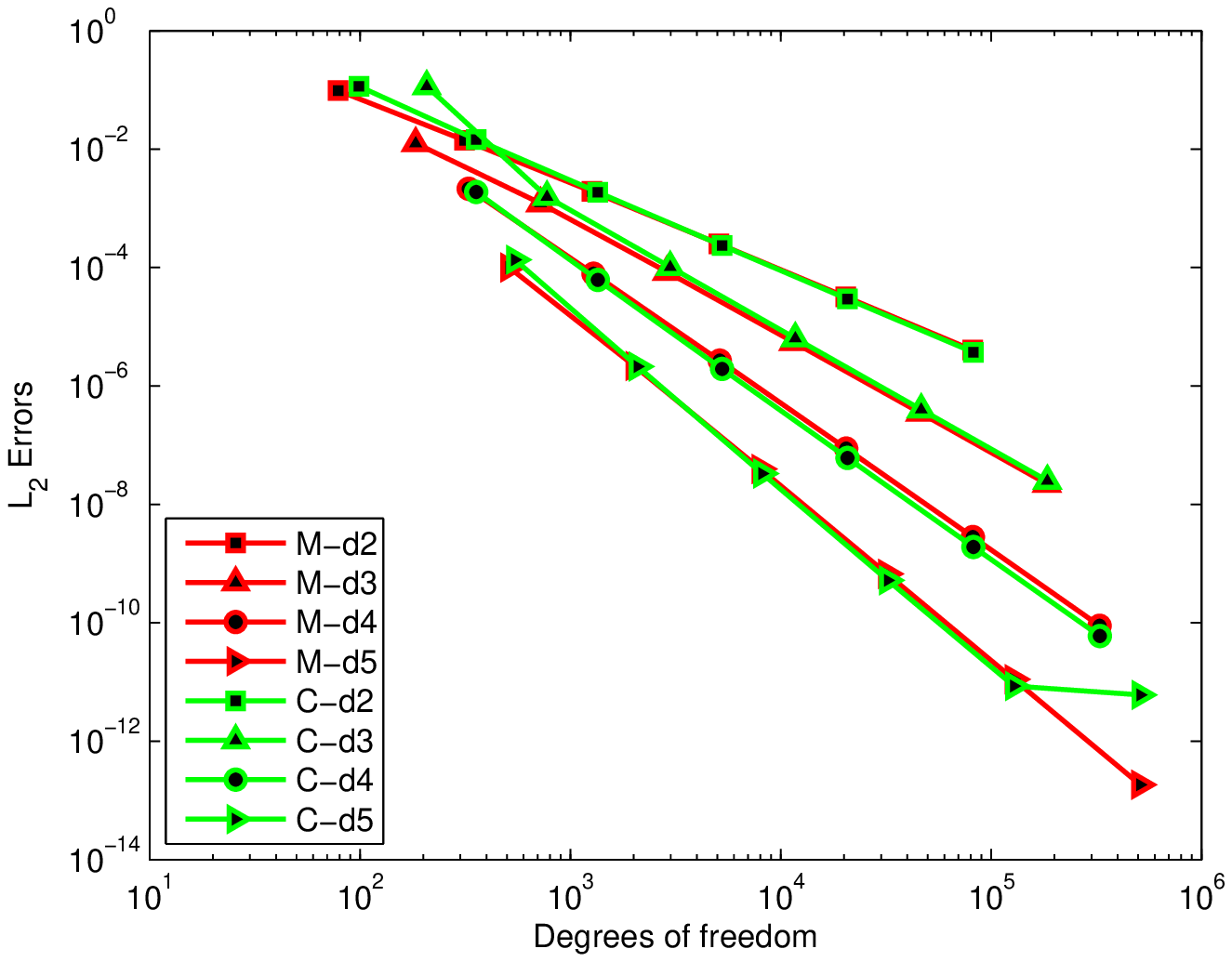} 
\end{tabular}
\caption{
$L_2$-errors in Example 3 using our method and COMSOL.
}\label{L2Test3}
\end{figure}

\begin{figure}
\centering
\begin{tabular}{c}
\includegraphics[height=0.6\textwidth]{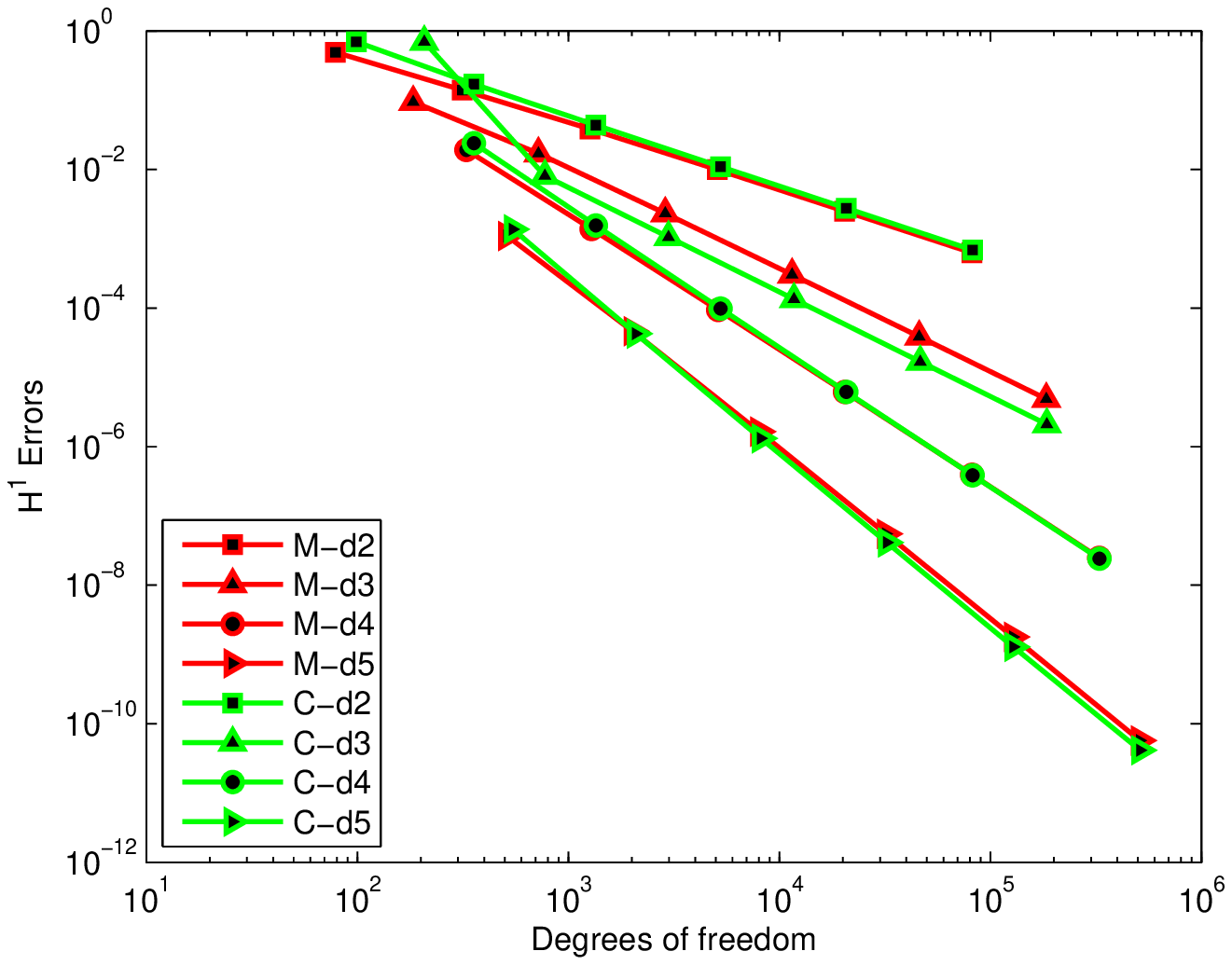} 
\end{tabular}
\caption{
$H^1$-errors in Example 3 using our method and COMSOL.
}\label{H1Test3}
\end{figure}

Similar to Example 1 we consider elements of degrees $d=2,3,4,5$ and compare our results with COMSOL. 
Figures~\ref{L2Test3} and \ref{H1Test3} depict the $L_2$ and $H^1$ errors for both methods. 
Again the numbers show the robust behavior of our method for different degrees while
confirming the error bounds of Theorem~\ref{bound1}. We do not consider the $p$-refinement for 
this example because the solution is a polynomial of degree 6 and hence lies in the approximation space
$S_{d,0}(\tri)$ for all $d\ge6$.

\end{document}